\documentclass[a4paper,11pt,reqno]{amsart}

\newtheorem{theorem}[equation]{Theorem}

\newtheorem{corollary}[equation]{Corollary}
\newtheorem{proposition}[equation]{Proposition}
\numberwithin{equation}{section}

\theoremstyle{definition}

\newtheorem{remark}[equation]{Remark}

\newcommand{\bN}{{\mathbb N}}
\newcommand{\bZ}{{\mathbb Z}}
\newcommand{\bR}{{\mathbb R}}

\newcommand{\bO}{{\mathbb O}}
\newcommand{\bP}{{\mathbb P}}

\newcommand{\frg}{{\mathfrak g}}
\newcommand{\frh}{{\mathfrak h}}
\newcommand{\frc}{{\mathfrak c}}
\newcommand{\frd}{{\mathfrak d}}
\newcommand{\frf}{{\mathfrak f}}
\newcommand{\fre}{{\mathfrak e}}
\newcommand{\frt}{{\mathfrak t}}
\newcommand{\frs}{{\mathfrak s}}
\newcommand{\fra}{{\mathfrak a}}
\newcommand{\frgo} {{\frg_{\bar 0}}}
\newcommand{\frguno} {{\frg_{\bar 1}}}
\newcommand{\frgoo}{{\frg_{(\bar 0,\bar 0)}}}
\newcommand{\frgouno}{{\frg_{(\bar 0,\bar 1)}}}
\newcommand{\frgunoo}{{\frg_{(\bar 1,\bar 0)}}}
\newcommand{\frgunouno}{{\frg_{(\bar 1,\bar 1)}}}
\newcommand{\calS}{{\mathcal S}}
\newcommand{\calSf}{{\mathcal S}_{\frf_4}}
\newcommand{\calSe}[1]{{\mathcal S}_{\fre_{#1}}}
\newcommand{\epsf}{\epsilon_{\frf_4}}
\newcommand{\epse}[1]{\epsilon_{\fre_{#1}}}
\newcommand{\calT}{{\mathcal T}}

\providecommand{\espan}[1]{\text{span}\left\{ #1\right\}}

 \DeclareMathOperator{\tri}{\mathfrak{tri}}
 \DeclareMathOperator{\fro}{\mathfrak{o}}
 \DeclareMathOperator{\frsl}{{\mathfrak{sl}}}
 \DeclareMathOperator{\frsp}{{\mathfrak{sp}}}
 \DeclareMathOperator{\frpsl}{{\mathfrak{psl}}}
 \DeclareMathOperator{\frpgl}{{\mathfrak{pgl}}}

 \DeclareMathOperator{\trace}{trace}
 
 \DeclareMathOperator{\Der}{Der}
 \DeclareMathOperator{\inder}{inder}
 \DeclareMathOperator{\End}{End}

\def\hregla{\hrule height.1pt}
\def\hreglabis{\hrule height .3pt depth -.2pt}
\def\hregleta{\hrule height .5pt}
\def\hreglon{\hrule height1pt}
\def\vreglon{\vrule height 12pt width1pt depth 4pt}
\def\vregleta{\vrule width .5pt}

\def\hreglafill{\leaders\hreglabis\hfill}
\def\hregletafill{\leaders\hregleta\hfill}
\def\vregla{\vrule width.1pt}

\begin{document}

\title{The Magic Square and Symmetric Compositions II}

\author{Alberto Elduque}
 \thanks{Supported by the Spanish Ministerio de Ciencia y
Tecnolog\'{\i}a and FEDER (BFM 2001-3239-C03-03), by the Ministerio de
Educaci\'on y Ciencia and FEDER (MTM 2004--08115-C04-02), and by the
Diputaci\'on General de Arag\'on (Grupo de Investigaci\'on de
\'Algebra)}
 \address{Departamento de Matem\'aticas, Universidad de
Zaragoza, 50009 Zaragoza, Spain}
 \email{elduque@unizar.es}




\begin{abstract}
The construction of Freudenthal's Magic Square, which contains the
exceptional simple Lie algebras of types $F_4,E_6,E_7$ and $E_8$, in
terms of symmetric composition algebras is further developed here.
The para-Hurwitz algebras, which form a subclass of the symmetric
composition algebras, will be defined, in the split case,  in terms
of the natural two dimensional module for the simple Lie algebra
$\frsl_2$. As a consequence, it will be shown how all the Lie
algebras in Freudenthal's Magic Square can be constructed, in a
unified way, using copies of $\frsl_2$ and of its natural module.
\end{abstract}

\maketitle


\section{Introduction}
The exceptional simple Lie algebras in Killing-Cartan's
classification are fundamental objects in many branches of
mathematics and physics. A lot of different constructions of these
objects have been given, many of which involve some nonassociative
algebras or triple systems.

In 1966 Tits gave a unified construction of the exceptional simple
Lie algebras which uses a couple of ingredients: a unital
composition algebra $C$ and a simple  Jordan algebra $J$ of degree
$3$ \cite{T}. At least in the split cases, this is a construction
which depends on two unital composition algebras, since the Jordan
algebra involved  consists of the $3\times 3$-hermitian matrices
over a unital composition algebra. Even though the construction is
not symmetric in the two composition algebras that are being used,
the outcome (the Magic Square) is indeed symmetric.

Over the years, more symmetric constructions have been given,
starting with a construction by Vinberg in 1966 \cite{OV}. Later
on, a quite general construction was given by Allison and Faulkner
\cite{AF} of Lie algebras out of structurable ones. In the
particular case of the tensor product of two unital composition
algebras, this construction provides another symmetric
construction of Freudenthal's Magic Square. Quite recently, Barton
and Sudbery \cite{BS1,BS2} (see also Landsberg and Manivel
\cite{LM1,LM2}) gave a simple recipe to obtain the Magic Square in
terms of two unital composition algebras and their triality Lie
algebras which, in perspective, is subsumed in Allison-Faulkner's
construction.

Let us recall that a \emph{composition algebra} is a triple
$(S,\cdot,q)$, where $(S,\cdot)$ is a (nonassociative) algebra
over a field $F$ with multiplication denoted by $x\cdot y$ for
$x,y\in S$, and where $q:S\rightarrow F$ is a regular quadratic
form (the \emph{norm}) such that, for any $x,y\in S$:
\begin{equation}\label{eq:comp}
 q(x\cdot y)=q(x)q(y).
\end{equation}

\emph{In what follows, the ground field $F$ will always be assumed
to be of characteristic $\ne 2$.}

Unital composition algebras (or Hurwitz algebras) form a
well-known class of algebras. Any Hurwitz algebra has finite
dimension equal to either $1$, $2$, $4$ or $8$. The
two-dimensional Hurwitz algebras are the quadratic \'etale
algebras over the ground field $F$, the four dimensional ones are
the generalized quaternion algebras, and the eight dimensional
Hurwitz algebras are called \emph{Cayley} algebras, and are
analogues to the classical algebra of octonions (for a nice survey
of the latter, see \cite{Baez}).

However, as shown in \cite{KMRT}, the triality phenomenon is
better dealt with by means of the so called \emph{symmetric
composition algebras}, instead of the classical unital composition
algebras. This led the author \cite{Ibero} to reinterpret the
Barton-Sudbery's construction in terms of two symmetric
composition algebras.

A composition algebra $(S,*,q)$ is said to be \emph{symmetric} if
it satisfies
\begin{equation}\label{eq:sym1}
 q(x*y,z)=q(x,y*z),
\end{equation}
where $q(x,y)=q(x+y)-q(x)-q(y)$ is the polar of $q$. In what
follows, any quadratic form and its polar will always be denoted
by the same letter. Equation \eqref{eq:sym1} is equivalent to
\begin{equation}\label{eq:sym2}
(x*y)*x=x*(y*x)=q(x)y
\end{equation}
for any $x,y\in S$. (See \cite[Ch.~VIII]{KMRT} for the basic facts
and notations.)

The classification of the symmetric composition algebras was
obtained in \cite{EM} (for characteristic $\ne 3$, see also
\cite[Ch.~VIII]{KMRT}) and in \cite{E1} (for characteristic $3$).

Given any  Hurwitz algebra $C$ with norm $q$, standard involution
$x\mapsto \bar x=q(x,1)1-x$, and multiplication denoted by
juxtaposition, the new algebra defined on $C$ but with
multiplication
$$
x\bullet y=\bar x\bar y,
$$
is a symmetric composition algebra, called the associated
\emph{para-Hurwitz algebra}. In dimension $1$, $2$ or $4$, any
symmetric composition algebra is a para-Hurwitz algebra, with a
few exceptions in dimension $2$ which are, nevertheless, forms of
para-Hurwitz algebras; while in dimension $8$, apart from the
para-Hurwitz algebras, there is a new family of symmetric
composition algebras termed \emph{Okubo algebras}.

If $(S,*,q)$ is any symmetric composition algebra, consider the
corresponding orthogonal Lie algebra
\[
\fro(S,q)=\{ d\in \End_F(S):
q\bigl(d(x),y\bigr)+q\bigl(x,d(y)\bigr)=0\ \forall x,y\in S\},
\]
and the subalgebra of $\fro(S,q)^3$ defined by
\begin{equation}\label{eq:triality}
 \tri(S,*,q)=\{(d_0,d_1,d_2)\in \fro(S,q)^3 :
 d_0(x*y)=d_1(x)*y+x*d_2(y)\ \forall x,y\in S\}.
\end{equation}

The map
\[
 \theta: \tri(S,*,q)\rightarrow \tri(S,*,q),\quad
 (d_0,d_1,d_2)\mapsto (d_2,d_0,d_1),
\]
is an automorphism of $\tri(S,*,q)$ of order $3$, the
\emph{triality automorphism}. Its fixed subalgebra is (isomorphic
to) the derivation algebra of the algebra $(S,*)$ which, if the
dimension is $8$ and the characteristic of the ground field is
$\ne 2,3$, is a simple Lie algebra of type $G_2$ in the
para-Hurwitz case and a simple Lie algebra of type $A_2$ (a form
of $\frsl_3$) in the Okubo case.

For any $x,y\in S$, the triple
\begin{equation}\label{eq:txy}
t_{x,y}=
 \left( \sigma_{x,y},\frac{1}{2}q(x,y)id-r_xl_y,
 \frac{1}{2}q(x,y)id-l_xr_y\right)
\end{equation}
is in $\tri(S,*,q)$, where $\sigma_{x,y}(z)=q(x,z)y-q(y,z)x$,
$r_x(z)=z*x$, and $l_x(z)=x*z$ for any $x,y,z\in S$.

\bigskip

The construction given in \cite{Ibero} starts with two symmetric
composition algebras $(S,*,q)$ and $(S',*,q')$. Then define
$\frg=\frg(S,S')$ to be the $\bZ_2\times\bZ_2$-graded
anticommutative algebra such that
$\frgoo=\tri(S,*,q)\oplus\tri(S',*,q')$,
 $\frgunoo=\frgouno=\frgunouno=S\otimes S'$.
(Unadorned tensor products are considered over the ground field
$F$.) For any $a\in S$ and $x\in S'$, denote by $\iota_i(a\otimes
x)$ the element $a\otimes x$ in $\frgunoo$ (respectively
$\frgouno$, $\frgunouno$) if $i=0$ (respectively, $i=1,2$). Thus
\begin{equation}\label{eq:gss'}
\frg=\frg(S,S')=\bigl(\tri(S,*,q)\oplus\tri(S',*,q')\bigr) \oplus
\bigl(\oplus_{i=0}^2\iota_i(S\otimes S')\bigr)\,.
\end{equation}
The anticommutative multiplication on $\frg$ is defined by means
of:
\begin{itemize}
\item $\frgoo$ is a Lie subalgebra of $\frg$,
\smallskip
\item $[(d_0,d_1,d_2),\iota_i(a\otimes
 x)]=\iota_i\bigl(d_i(a)\otimes x\bigr)$,
 $[(d_0',d_1',d_2'),\iota_i(a\otimes
 x)]=\iota_i\bigl(a\otimes d_i'(x)\bigr)$, for any $(d_0,d_1,d_2)\in
 \tri(S,*,q)$, $(d_0',d_1',d_2')\in \tri(S',*,q')$, $a\in S$ and
 $x\in S'$.
\smallskip
\item $[\iota_i(a\otimes x),\iota_{i+1}(b\otimes y)]=
 \iota_{i+2}\bigl((a*b)\otimes (x*y)\bigr)$ (indices modulo $3$), for any
 $a,b\in S$, $x,y\in S'$.
\smallskip
\item $[\iota_i(a\otimes x),\iota_i(b\otimes y)]=
 q'(x,y)\theta^i(t_{a,b})+q(a,b)\theta'^i(t'_{x,y})$, for any
 $i=0,1,2$, $a,b\in S$ and $x,y\in S'$, where $t_{a,b}\in
 \tri(S,*,q)$ (respectively $t'_{x,y}\in \tri(S',*,q')$) is the
 element in \eqref{eq:txy} for $a,b\in S$ (resp. $x,y\in S'$) and $\theta$
 (resp. $\theta'$) is the triality automorphism of $\tri(S,*,q)$ (resp.
 $\tri(S',*,q')$).

\end{itemize}

\bigskip

The main result in \cite{Ibero} asserts that, with this
multiplication, $\frg(S,S')$ is a Lie algebra and, if the
characteristic of the ground field is $\ne 2,3$, Freudenthal's
Magic Square is recovered (Table \ref{ta:magic}).

\begin{table}[h!]
$$ \vbox{\offinterlineskip
 \halign{$#$&\hfil\quad$#$\quad\hfil&%
 \vrule height 12pt depth 4pt #%
 &\hfil\quad$#$\quad\hfil&\hfil\quad$#$\quad\hfil
 &\hfil\quad$#$\quad\hfil&\hfil\quad$#$\ \hfil\cr
 &&\omit&\multispan4{\hfil$\dim S$\hfil}\cr
 \noalign{\smallskip}
 \strut &&width 0pt&1&2&4&8\cr
 &&\multispan5{\hrulefill}\cr
 &1&&A_1&A_2&C_3&F_4\cr
 \strut&2&&A_2&\omit$A_2\oplus A_2$&A_5&E_6\cr
 \strut\smash{\raise 6pt\hbox{$\dim S'$}}&4&&C_3&A_5&D_6&E_7\cr
 \strut&8&&F_4&E_6&E_7&E_8\cr}}
$$
\smallskip
\caption{The Magic Square}\label{ta:magic}
\end{table}

\medskip

In \cite{Mejico} it is proved that if $(S,*,q)$ is an Okubo
algebra with an idempotent (this is always the case if the ground
field does not admit cubic extensions), then there is a
para-Cayley algebra $(\hat S,\bullet,q)$ defined on the same
vector space $S$ such that $\frg(S,S')$ is isomorphic to
$\frg(\hat S,S')$, for any $S'$.

\bigskip

The purpose of this paper is to delve in the construction
$\frg(S,S')$ given in \cite{Ibero}, for two split para-Hurwitz
algebras. It turns out that the split para-Hurwitz algebras of
dimension $4$ and $8$ (para-quaternion and para-octonion algebras)
can be constructed using a very simple ingredient: the two
dimensional module for the three dimensional simple split Lie
algebra $\frsl_2$. This will be shown in Section 2.

As a consequence, all the split Lie algebras in the Magic Square
will be constructed in very simple terms using copies of $\frsl_2$
and of its natural module. The general form will be
\[
\frg=\oplus_{\sigma\in\calS}V(\sigma)\,,
\]
where $\calS$ will be a set of subsets of $\{1,\ldots,n\}$, for
some $n$, $V(\emptyset)=\oplus_{i=1}^n\frsl_2$, a direct sum of
copies of $\frsl_2$, and for
$\emptyset\ne\sigma=\{i_1,\ldots,i_r\}\in\calS$,
$V(\sigma)=V_{i_1}\otimes\cdots\otimes V_{i_r}$, with $V_i$ the
two-dimensional natural module for the $i^\text{th}$ copy of
$\frsl_2$ in $V(\emptyset)$. The Lie bracket will appear in terms
of natural contractions. The precise formulas will be given in
Section 3. Surprisingly, the real division algebra of octonions
will appear in these constructions, but in an unexpected way.

Section 4 will be devoted to show how the models obtained in
Section 3 of the exceptional simple split Lie algebras give models
too of the exceptional simple split Freudenthal triple systems
which, together with simple Jordan algebras, were the basic tools
used by Freudenthal to construct the Lie algebras in the Magic
Square.

\section{Split para-quaternions and para-octonions}

Let $V$ be a two dimensional vector space over a ground field $F$
(of characteristic $\ne 2$), endowed with a nonzero skew-symmetric
bilinear form $\langle .\vert .\rangle$. Consider the symplectic
Lie algebra
\begin{equation}\label{eq:gammaab}
\frsp(V)=\espan{\gamma_{a,b}=\langle a\vert .\rangle b + \langle
b\vert .\rangle a : a,b\in V}\,,
\end{equation}
which coincides with $\frsl(V)$ (endomorphisms of zero trace). The
bilinear form allows $Q=V\otimes V$ to be identified with
$\End_F(V)$ (the split quaternion algebra over $F$) by means of:
\[
\begin{split}
V\otimes V&\longrightarrow \End_F(V)\\
a\otimes b\,&\mapsto\ \langle a\vert .\rangle b: v\mapsto \langle
a\vert v\rangle b\, .
\end{split}
\]
Since
\[
\langle a\otimes b(v)\vert w\rangle =\langle a\vert v\rangle
  \langle b\vert w\rangle =-\langle v\vert b\otimes a(w)\rangle\,,
\]
it follows that for any $a,b\in V$,
\[
\overline{a\otimes b}=-b\otimes a\,,
\]
where $f\mapsto \bar f$ is the symplectic involution of
$\End_F(V)$ relative to $\langle .\vert .\rangle$ (the standard
involution as a quaternion algebra). Also, the norm as a
quaternion algebra is given by $q(x)=\det x=\frac{1}{2}\trace
x\bar x$; so its polar form becomes
\[
\begin{split}
q(a\otimes b,c\otimes d)
 &=\trace\bigl((a\otimes b)(\overline{c\otimes d})\bigr)\\
 &=-\trace\bigl((a\otimes b)(d\otimes c)\bigr)\\
 &=-\langle a\vert c\rangle\langle d\vert b\rangle\\
 &=\langle a\vert c\rangle\langle b\vert d\rangle\, .
 \end{split}
\]
(The natural symmetric bilinear form on $V\otimes V$ induced from
$\langle .\vert .\rangle$.)

The split Cayley algebra is obtained by means of the
Cayley-Dickson doubling process (see \cite[Section 2]{Jacobson}):
$C=Q\oplus Q$ with multiplication, standard involution and norm
given by:
\[
\left\{%
\begin{aligned}
&(x,y)(x',y')=(xx'-\overline{y'}y, y'x+y\overline{x'}),\\
&\overline{(x,y)}=(\bar x,-y),\\
&q\bigl((x,y)\bigr)=q(x)+q(y),
\end{aligned}
\right.
\]
for any $x,y,x',y'\in Q$. Therefore, $C=V\otimes
V\,\oplus\,V\otimes V$ with
\[
\left\{%
\begin{aligned}
&(a\otimes b,c\otimes d)(a'\otimes b',c'\otimes d')\\
&\qquad =\Bigl(\langle a\vert b'\rangle a'\otimes b-
     \langle d\vert d'\rangle c\otimes c',
     -\langle b\vert c'\rangle a\otimes d'-
     \langle c\vert a'\rangle b'\otimes d\Bigr),\\[3pt]
&\overline{(a\otimes b,c\otimes d)}=-(b\otimes a,c\otimes d),\\
&q\Bigl((a\otimes b,c\otimes d),(a'\otimes b',c'\otimes d')\Bigr)
     =\langle a\vert a'\rangle\langle b\vert b'\rangle +
      \langle c\vert c'\rangle\langle d\vert d'\rangle\, ,
\end{aligned}
\right.
\]
for any $a,b,c,d,a',b',c',d'\in V$.

The multiplication $x\bullet y=\bar x\bar y$ in the associated
para-Hurwitz algebra (the split para-octonions) $S_8$ takes the
form:
\begin{equation}\label{eq:bullet}
\begin{split}
(a\otimes b,c&\otimes d)\bullet (a'\otimes b',c'\otimes d')\\
  &=(b\otimes a,c\otimes d)(b'\otimes a',c'\otimes d')\\
  &=\Bigl(\langle b\vert a'\rangle b'\otimes a-
     \langle d\vert d'\rangle c\otimes c',
     -\langle a\vert c'\rangle b\otimes d'-
     \langle c\vert b'\rangle a'\otimes d\Bigr).
\end{split}
\end{equation}
The split para-quaternion algebra $S_4$ is just the subalgebra
consisting of the first copy of $V\otimes V$ in $S_8$.

The Lie algebra $\frsp(V)^4$ acts naturally on $S_8$, where the
$i^\text{th}$ component of $\frsp(V)^4$ acts on the $i^\text{th}$
copy of $V$ in $S_8=V\otimes V\,\oplus\,V\otimes V$. This gives an
embedding into the orthogonal Lie algebra of $S_8$ relative to
$q$:
\begin{equation}\label{eq:rhospv4}
\rho:\frsp(V)^4\longrightarrow \fro(S_8,q).
\end{equation}
Actually, this is an isomorphism of $\frsp(V)^4$ onto the
subalgebra $\fro(V\otimes V)\oplus\fro(V\otimes V)$ of
$\fro(S_8,q)$, which is the even part of $\fro(S_8,q)$ relative to
the $\bZ_2$-grading given by the orthogonal decomposition
$S_8=V\otimes V\,\perp\,V\otimes V$.

Consider also the linear map (denoted by $\rho$ too):
\begin{equation}\label{eq:rhov4}
\rho :\ V^{\otimes 4}\longrightarrow \fro(S_8,q)\,,
\end{equation}
such that
\begin{equation}\label{eq:rhoviwi}
\begin{split}
\rho(v_1\otimes v_2\otimes v_3\otimes v_4)&
 \bigl((w_1\otimes w_2,w_3\otimes w_4)\bigr)\\
 &=
 \Bigl( \langle v_3\vert w_3\rangle\langle v_4\vert w_4\rangle
    v_1\otimes v_2,
 -\langle v_1\vert w_1\rangle\langle v_2\vert w_2\rangle
    v_3\otimes v_4\Bigr)\,,
\end{split}
\end{equation}
for any $v_i,w_i\in V$, $i=1,2,3,4$. Consider for any $x,y\in S_8$
the linear map
\[
\sigma_{x,y}=q(x,.)y-q(y,.)x
\]
(these maps span $\fro(S_8,q)$). Then,
\begin{equation}\label{eq:rhosigma}
\rho(v_1\otimes v_2\otimes v_3\otimes v_4)=
 -\sigma_{(v_1\otimes v_2,0),(0,v_3\otimes v_4)}\, .
\end{equation}
Observe that $\rho(v_1\otimes v_2\otimes v_3\otimes v_4)$ fills
the odd part of the $\bZ_2$-grading of $\fro(S_8,q)$ mentioned
above.

Now, for any $v_i,w_i\in V$, $i=1,2,3,4$; a straightforward
computation, using that $\langle a\vert b\rangle c+\langle b\vert
c\rangle a+\langle c\vert a\rangle b=0$ for any $a,b,c\in V$
(since $\dim V=2$ and the expression is skew symmetric on its
arguments), gives:
\begin{equation}\label{eq:[rho,rho]}
\begin{split}
\Bigl[ \rho(v_1\otimes v_2\otimes v_3\otimes v_4)&,
  \rho(w_1\otimes w_2\otimes w_3\otimes w_4)\Bigr]\\
&=\frac{1}{2} \sum_{i=1}^4
 \Bigl(\prod_{j\ne i}\langle v_j\vert w_j\rangle\Bigr)
 \nu_i(\gamma_{v_i,w_i})\, \in\rho\bigl(\frsp(V)^4\bigr)\,,
\end{split}
\end{equation}
where $\nu_i:\frsp(V)\hookrightarrow \frsp(V)^4$ denotes the
inclusion into the $i^\text{th}$ component.

In the same vein, for any $v_i,w_i\in V$, $i=1,2,3,4$.
\begin{equation}\label{eq:rhosigma2}
\left\{%
\begin{aligned}
\sigma_{(v_1\otimes v_2,0),(w_1\otimes w_2,0)}&
 =\frac{1}{2}\Bigl(\langle v_2\vert w_2\rangle
 \nu_1(\gamma_{v_1,w_1}) +\langle v_1\vert w_1\rangle
  \nu_2(\gamma_{v_2,w_2})\Bigr),\\
\sigma_{(0,v_3\otimes v_4),(0,w_3\otimes w_4)}&
 =\frac{1}{2}\Bigl(\langle v_4\vert w_4\rangle
 \nu_3(\gamma_{v_3,w_3}) +\langle v_3\vert w_3\rangle
  \nu_4(\gamma_{v_4,w_4})\Bigr).
\end{aligned}\right.
\end{equation}

The following result (see also \cite[3.4]{LM2}) summarizes most of
the above arguments.

\begin{proposition}\label{pr:d4}
The vector space
\begin{equation}\label{eq:d4}
\frd_4=\frsp(V)^4\oplus V^{\otimes 4}\,,
\end{equation}
with anticommutative multiplication given by:

\begin{itemize}

\item $\frsp(V)^4$ is a Lie subalgebra of $\frd_4$,
\smallskip
\item for any  $s_i\in
 \frsp(V)$ and $v_i\in V$, $i=1,2,3,4$,
 \[
 \begin{split}
 [(s_1,s_2,s_3,s_4)&,v_1\otimes v_2\otimes v_3\otimes v_4]\\
 &=
 s_1(v_1)\otimes v_2\otimes v_3\otimes v_4 +
 v_1\otimes s_2(v_2)\otimes v_3\otimes v_4 \\
 &\qquad+
 v_1\otimes v_2\otimes s_3(v_3)\otimes v_4 +
 v_1\otimes v_2\otimes v_3\otimes s_4(v_4),
 \end{split}
 \]
\smallskip
\item for any $v_i,w_i\in V$, $i=1,2,3,4$,
 \[
 [v_1\otimes v_2\otimes v_3\otimes v_4,
         w_1\otimes w_2\otimes w_3\otimes w_4]=
 \frac{1}{2}\sum_{i=1}^4
 \bigl(\prod_{j\ne i}\langle v_j\vert w_j\rangle\bigr)
 \nu_i(\gamma_{v_i,w_i})\,,
 \]
\end{itemize}
is a Lie algebra isomorphic to $\fro(S_8,q)$ by means of the
linear map $\rho$ defined by \eqref{eq:rhospv4} and
\eqref{eq:rhov4}.
\end{proposition}

Note that $\frd_4$ is generated by $V^{\otimes 4}$ and that the
decomposition $\frd_4=\frsp(V)^4\oplus V^{\otimes 4}$ is a
$\bZ_2$-grading.

Let us consider now, as in \cite[3.4]{LM2}, the order three
automorphism $\theta:\frd_4\rightarrow \frd_4$ such that:
\begin{equation}\label{eq:theta}
\left\{%
\begin{aligned}
&\theta\bigl((s_1,s_2,s_3,s_4)\bigr)=(s_3,s_1,s_2,s_4),\\
&\theta(v_1\otimes v_2\otimes v_3\otimes v_4)=
  v_3\otimes v_1\otimes v_2\otimes v_4,
\end{aligned}
\right.
\end{equation}
for any $s_i\in \frsp(V)$ and $v_i\in v$, $i=1,2,3,4$. Then:

\begin{proposition}[Local triality] For any $f\in \frd_4$ and
any $x,y\in S_8$:
\[
\rho(f)(x\bullet y)=\Bigl(\rho(\theta^{-1}(f))(x)\Bigr)\bullet y
  +x\bullet\Bigl(\rho(\theta^{-2}(f))(y)\Bigr)\,.
\]
\end{proposition}
\begin{proof}
It is enough to prove this for generators of $\frd_4$ and of
$S_8$, and hence for $f=v_1\otimes v_2\otimes v_3\otimes v_4$,
$x=(a_1\otimes a_2,a_3\otimes a_4)$, $y=(b_1\otimes b_2,b_3\otimes
b_4)$, with $v_i,a_i,b_i\in V$, $i=1,2,3,4$, and this is a
straightforward computation.
\end{proof}

The triality Lie algebra of $(S_8,\bullet,q)$ is
\[
\begin{split}
&\tri(S_8,\bullet,q)\\
&\quad=\{ (g_0,g_1,g_2)\in \fro(S_8,q) : g_0(x\bullet
y)=g_1(x)\bullet y+x\bullet g_2(y),\ \forall x,y\in S_8\}\, ,
\end{split}
\]
and the projection on any of its components gives an isomorphism
among the Lie algebras $\tri(S_8,\bullet,q)$ and $\fro(S_8,q)$.
Let $\rho_i=\rho\circ\theta^{-i}$, $i=0,1,2$. By dimension count,
the previous Proposition immediately implies:

\begin{corollary}\label{co:trirho}
$\tri(S_8,\bullet,q)=\left\{\bigl((\rho_0(f),\rho_1(f),\rho_2(f)\bigr):
f\in \frd_4\right\}$.
\end{corollary}

Denote by $V_i$ the $\frsp(V)^4$-module $V$ on which only  the
$i^\text{th}$ component acts: $(s_1,s_2,s_3,s_4).v=s_i(v)$ for any
$s_i\in\frsp(V)$, $i=1,2,3,4$, and $v\in V$. Also, denote by
$\iota_i(S_8)$ the $\frd_4$-module associated to the
representation $\rho_i:\frd_4\rightarrow \fro(S_8,q)$. Then, as
modules for $\frsp(V)^4$:
\begin{equation}\label{eq:iota012}
\begin{split}
\iota_0(S_8)&=V_1\otimes V_2\,\oplus\,V_3\otimes V_4,\\
\iota_1(S_8)&=V_2\otimes V_3\,\oplus\,V_1\otimes V_4,\\
\iota_2(S_8)&=V_3\otimes V_1\,\oplus\,V_2\otimes V_4\,(\simeq
V_1\otimes V_3\,\oplus\,V_2\otimes V_4).
\end{split}
\end{equation}

\begin{remark} For the split para-quaternion algebra $S_4=V\otimes
V$, by restriction we obtain:
\[
\tri(S_4,\bullet,q)=\left\{
\bigl(\tilde\rho_0(f),\tilde\rho_1(f),\tilde\rho_2(f)\bigr) :
f\in\frsp(V)^3\right\},
\]
where $\tilde\rho_i$ is obtained by restriction of $\rho_i$:
\[
\tilde\rho_i\bigl((s_1,s_2,s_3)\bigr)
 =\rho_i\bigl((s_1,s_2,s_3,0)\bigr)\vert_{V\otimes V}\,.
\]
\end{remark}

\begin{remark}
The construction of $\frd_4$ in Proposition \ref{pr:d4} makes it
clear an action of the symmetric group $S_3$ on $\frd_4$ leaving
fixed the last copy of $\frsp(V)$ and of $V$, as shown by the
action of $\theta$ in \eqref{eq:theta}. As such, $\frd_4$ is the
natural example of a \emph{Lie algebra with triality}, as defined
in \cite{Grishkovtriality}.

Moreover,  the fixed subalgebra by $\theta$ in \eqref{eq:theta} is
a direct sum $\bigl(\frsp(V)\oplus\frsp(V)\bigr)\oplus S^3V\otimes
V$, where the first copy of $\frsp(V)$ is the diagonal subalgebra
in the direct sum of the first three copies of $\frsp(V)$ in
\eqref{eq:d4}, $S^3(V)$ is the module of symmetric tensors in
$V^{\otimes 3}$ (the tensor product of the first three copies of
$V$ in \eqref{eq:d4}). If the characteristic is $\ne 2,3$, this is
the split simple Lie algebra of type $G_2$.
\end{remark}

\bigskip

The following notation will be useful in the sequel. For any $n\in
\bN$ and any subset $\sigma\subseteq \{1,2,\ldots,n\}$ consider
the $\frsp(V)^n$-modules given by
\[
V(\sigma)=\begin{cases}%
\frsp(V)^n& \text{if $\sigma=\emptyset$,}\\
V_{i_1}\otimes \cdots\otimes V_{i_r}&\text{if
$\sigma=\{i_1,\ldots,i_r\}$, $1\leq i_1<\cdots<i_r\leq n$.}
\end{cases}
\]
As before, $V_i$ denotes the module $V$ for the $i^\text{th}$
component of $\frsp(V)^n$, annihilated by the other $n-1$
components.

Identify any subset $\sigma\subseteq \{1,\ldots,n\}$ with the
element $(\sigma_1,\ldots,\sigma_n)\in \bZ_2^n$ such that
$\sigma_i=1$ if and only if $i\in \sigma$. Then for any
$\sigma,\tau\in\bZ_2^n$, consider the natural
$\frsp(V)^n$-invariant maps
\begin{equation}\label{eq:phisigmatau}
\varphi_{\sigma,\tau}:V(\sigma)\times V(\tau)\longrightarrow
V(\sigma+\tau)
\end{equation}
defined as follows:

\begin{itemize}

\item
If $\sigma\ne \tau$ and $\sigma\ne\emptyset\ne\tau$, then
$\varphi_{\sigma,\tau}$ is obtained by contraction, by means of
$\langle .\vert .\rangle$ in the indices $i$ with
$\sigma_i=1=\tau_i$. Thus, for instance,
\[
\varphi_{\{1,2,3\},\{1,3,4\}}
 (v_1\otimes v_2\otimes v_3,w_1\otimes w_3\otimes w_4)=
 \langle v_1\vert w_1\rangle\langle v_3\vert w_3\rangle
     v_2\otimes w_4
\]
for any $v_1,w_1\in V_1$, $v_2\in V_2$, $v_3,w_3\in V_3$ and
$w_4\in V$.
\smallskip
\item
$\varphi_{\emptyset,\emptyset}$ is the Lie bracket in
$\frsp(V)^n$.
\smallskip
\item
For any $\sigma\ne\emptyset$, $\varphi_{\emptyset,\sigma} =
-\varphi_{\sigma,\emptyset}$ is given by the natural action of
$\frsp(V)^n$ on $V(\sigma)$. Thus, for instance,
\[
\varphi_{\emptyset,\{1,3\}}\bigl((s_1,\ldots,s_n),
 v_1\otimes v_3\bigr)=s_1(v_1)\otimes v_3 +
   v_1\otimes s_3(v_3),
\]
for any $s_i\in \frsp(V)$, $i=1,\ldots,n$, and $v_1\in V_1$,
$v_3\in V$.
\smallskip
\item
Finally, for any $\sigma\ne \emptyset$, $\varphi_{\sigma,\sigma}$
is given altering slightly  \eqref{eq:[rho,rho]}:
\[
\varphi_{\sigma,\sigma} (v_{i_1}\otimes \cdots\otimes v_{i_r},
   w_{i_1}\otimes \cdots \otimes w_{i_r})=
   -\frac{1}{2}\sum_{j=1}^r\Bigl(
   \prod_{k\ne j}\langle v_{i_k}\vert w_{i_k}\rangle
   \nu_{i_j}(\gamma_{v_{i_j},w_{i_j}})\Bigr)\,.
\]
(The minus sign here is useful in getting nicer formulae later
on.)
\end{itemize}

\bigskip
With this notation and for $n=4$,
\begin{equation}\label{eq:d4bis}
\frd_4=V(\emptyset)\oplus V(\{1,2,3,4\})\,,
\end{equation}
with
\begin{equation}\label{eq:xsigmaytaud4}
[x_\sigma,y_\tau]=\epsilon_{\frd_4}(\sigma,\tau)\varphi_{\sigma,\tau}(x_\sigma,y_\tau)
\end{equation}
for any $\sigma,\tau\in\Bigl\{\emptyset,\{1,2,3,4\}\Bigr\}$ where,
for $\sigma=\{ 1,2,3,4\}$,
\[
\epsilon_{\frd_4}(\emptyset,\emptyset)=\epsilon_{\frd_4}(\emptyset,\sigma)
=\epsilon_{\frd_4}(\sigma,\emptyset)=1\quad\text{and}\quad
\epsilon_{\frd_4}(\sigma,\sigma)=-1.
\]

\section{Exceptional Lie algebras}
In this section, the construction of the split para-quaternions
and para-octonions given in Section 2, together with the
description given there of $\frd_4$ and of $\tri(S_8,\bullet,q)$,
will be used to provide constructions of the exceptional simple
Lie algebras, which depend only on copies of $\frsp(V)$ and of
$V$. (Notations as in the previous section.)

\medskip

\subsection{$\mathbf{F_4}$} \quad\null

 Let us start with the split Lie algebra of type $F_4$.
Because of \eqref{eq:gss'}, \eqref{eq:iota012} and \eqref{eq:d4},
and identifying $\tri(S_8,,\bullet,q)$ with $\frd_4$ by means of
$\rho_0$ (Corollary \ref{co:trirho}), one has:
\begin{equation}\label{eq:f4}
\frf_4=\frg(S_8,F)=\oplus_{\sigma\in\calS_{\frf_4}}V(\sigma)\,,
\end{equation}
where
\[
\calSf=\Bigl\{\emptyset,\{1,2,3,4\},\{1,2\},\{2,3\},\{1,3\},
 \{3,4\},\{1,4\},\{2,4\}\Bigr\}\subseteq 2^{\{1,2,3,4\}}.
\]
By $\frsp(V)^4$-invariance of the Lie bracket in $\frf_4$, it
follows that
\begin{equation}\label{eq:xsigmaytau4}
[x_\sigma,y_\tau]=\epsf(\sigma,\tau)
\varphi_{\sigma,\tau}(x_\sigma,y_\tau)\,,
\end{equation}
for any $\sigma,\tau\in \calSf$, $x_\sigma\in V(\sigma)$ and
$y_\tau\in V(\tau)$; for a suitable map
$\epsf:\calSf\times\calSf\rightarrow F$.

Actually, the product $\bullet$ in \eqref{eq:bullet} becomes in
$\frg(S_8,F)$ a bilinear map
\[
\iota_0(S_8\otimes F)\times \iota_1(S_8\otimes F)\rightarrow
\iota_2(S_8\otimes F)
\]
(and cyclically); that is, a bilinear map
\[
\bigl(V_1\otimes V_2\,\oplus\,V_3\otimes V_4\bigr)\times
\bigl(V_2\otimes V_3\,\oplus\,V_1\otimes V_4\bigr)\rightarrow
\bigl(V_3\otimes V_1\,\oplus\,V_2\otimes V_4\bigr)
\]
given by
\[
\begin{split}
\Bigl((a_1&\otimes a_2,a_3\otimes a_4),(b_2\otimes b_3,b_1\otimes
b_4)\Bigr)\\
&\mapsto \Bigl(\langle a_2\vert b_2\rangle b_3\otimes a_1
 -\langle a_4\vert b_4\rangle a_3\otimes b_1,
 -\langle a_1\vert b_1\rangle a_2\otimes b_4 -
 \langle a_3\vert b_3\rangle b_1\otimes a_4\Bigr)\,,
\end{split}
\]
and cyclically. Thus it consists of maps $\varphi_{\sigma,\tau}$
scaled by $\pm 1$. Also, equation \eqref{eq:rhoviwi} determines
the Lie bracket in $\frg(S_8,F)$ between $V(\{1,2,3,4\})$ and
$V(\{1,2\})\oplus V(\{3,4\})$. By cyclic symmetry using $\rho_1$
and $\rho_2$ one gets all the brackets between $V(\{1,2,3,4\})$
and the $V(\sigma)$'s with $\sigma\in\calSf$ consisting of two
elements. Finally, for any $x,y\in S_8$, $[\iota_0(x\otimes
1),\iota_0(y\otimes 1)]=2t_{x,y}$ \eqref{eq:txy}, since
$q'(1,1)=2$ in the symmetric composition algebra $F$. But under
$\rho_0=\rho$, $t_{x,y}$ corresponds to
$\sigma_{x,y}\in\fro(S_8,q)$. Then equations \eqref{eq:rhosigma}
and \eqref{eq:rhosigma2} determine the corresponding values of
$\epsf$. By cyclic symmetry, one completes the information about
$\epsf$, which is displayed on Table \ref{ta:f4}.

\begin{table}[h!]
$$
\vbox{\offinterlineskip \halign{\hfil$#$\hfil\enspace\vreglon
 &\enspace\hfil$#$\hfil\enspace\vregla
 &\enspace\hfil$#$\hfil\enspace\vregla
 &\enspace\hfil$#$\hfil\enspace\vregla
 &\enspace\hfil$#$\hfil\enspace&\hskip .2pt\vregla#
 &&\enspace\hfil$#$\hfil\enspace\vregla\cr
 \omit\enspace\hfil\raise2pt\hbox{$\epsf$}%
  \enspace\hfil\vrule width1pt depth 4pt
   &\enspace\emptyset\enspace&\{1,2\}&\{2,3\}&\{1,3\}&\omit\vregleta
     &\omit$\{1,2,3,4\}$\vregla&\{3,4\}&\{1,4\}&\{2,4\}\cr
 \noalign{\hreglon}
 \emptyset&1&1&1&1&\omit\vregleta&1&1&1&1\cr
 \noalign{\hregla}
 \{1,2\}&1&-2&1&1&\omit\vregleta&1&-2&-1&-1\cr\noalign{\hregla}
 \{2,3\}&1&1&-2&1&\omit\vregleta&1&-1&-2&-1\cr\noalign{\hregla}
 \{1,3\}&1&1&1&-2&\omit\vregleta&1&-1&-1&-2\cr
 \multispan6{\hregletafill}&\multispan4{\hreglafill}\cr
 \{1,2,3,4\}&1&-1&-1&-1&&-1&1&1&1\cr\noalign{\hregla}
 \{3,4\}&1&2&-1&-1&&-1&-2&-1&-1\cr\noalign{\hregla}
 \{1,4\}&1&-1&2&-1&&-1&-1&-2&-1\cr\noalign{\hregla}
 \{2,4\}&1&-1&-1&2&&-1&-1&-1&-2\cr
 \noalign{\hregla}
}}
$$
\caption{$\epsf$}\label{ta:f4}
\end{table}

That is, $\epsf$ takes the following values for
$\sigma,\tau\subseteq \{1,2,3,4\}$:
\begin{equation}\label{eq:eps4}
\epsf(\sigma,\tau)=
\begin{cases}
\phantom{-}1&\text{if $\sigma$ or $\tau$ is empty,}\\[5pt]
-1&\text{if $\sigma=\tau=\{1,2,3,4\}$,}\\
-2&\text{if $\emptyset\ne\sigma=\tau\ne \{1,2,3,4\}$,}\\[5pt]
\ 2&\text{if $4\in\sigma= \{1,2,3,4\}\setminus \tau$,
           $\tau\ne\emptyset$,}\\
-2&\text{if $4\in\tau= \{1,2,3,4\}\setminus \sigma$,
           $\sigma\ne\emptyset$,}\\[5pt]
\ 1&\text{if $4\not\in \sigma\cup\tau$,}\\
-1&\text{if $4\in\sigma\cup\tau$, $\sigma\ne\emptyset\ne\tau$.}
\end{cases}
\end{equation}

\bigskip

The split Lie algebra $C_3$ is realized as $\frc_3=\frg(S_4,F)$, a
subalgebra of  $\frf_4=\frg(S_8,F)$ (see \cite[Theorem
2.6]{Mejico}). Thus, by restriction, one gets
\begin{equation}\label{eq:c3}
\begin{split}
\frc_3&=\frsp(V)^3\oplus\, V_1\otimes V_2\,\oplus\, V_2\otimes V_3\,\oplus
    \, V_1\otimes V_3 \\
  &=\oplus_{\sigma\in\calS_{\frc_3}}V(\sigma)\,,
\end{split}
\end{equation}
where
\[
\calS_{\frc_3}=\Bigl\{\emptyset,\{1,2\},\{2,3\},\{1,3\}\Bigr\}\subseteq
2^{\{1,2,3\}}.
\]
and the multiplication is given by
\begin{equation}\label{eq:xsigmaytauC3}
[x_\sigma,y_\tau]=\epsilon_{\frc_3}(\sigma,\tau)
\varphi_{\sigma,\tau}(x_\sigma,y_\tau)\,,
\end{equation}
for any $\sigma,\tau\in \calS_{\frc_3}$, $x_\sigma\in V(\sigma)$
and $y_\tau\in V(\tau)$; where
$\epsilon_{\frc_3}:\calS_{\frc_3}\times\calS_{\frc_3}\rightarrow
F$ is given by the left top corner of Table \ref{ta:f4}. That is,
\[
\epsilon_{\frc_3}(\sigma,\tau)=
\begin{cases}
-2&\text{if $\sigma=\tau\ne\emptyset$,}\\
1&\text{otherwise.}
\end{cases}
\]

\medskip

\subsection{$\mathbf{E_7}$} \quad\null

The split simple Lie algebra $E_7$ is realized as
$\fre_7=\frg(S_8,S_4)$. Here there are $4$ copies of $V$ involved
in $S_8$ and another $3$ copies for $S_4$ \eqref{eq:c3}. The
indices $1,2,3,4$ will be used for $S_8$ and the indices $5,6,7$
for $S_4$. Therefore, with the same arguments as above,
\begin{equation}\label{eq:e7}
\fre_7=\oplus_{\sigma\in\calSe{7}}V(\sigma)\,,
\end{equation}
where
\[
\begin{split}
\calSe{7}&=\Bigl\{\emptyset,\{1,2,5,6\},\{2,3,6,7\},\{1,3,5,7\},\\
  &\qquad \{1,2,3,4\},\{3,4,5,6\},\{1,4,6,7\},\{2,4,5,7\}\Bigr\}\\
  &\subseteq 2^{\{1,2,3,4,5,6,7\}}.
\end{split}
\]
and the multiplication is given by
\begin{equation}\label{eq:xsigmaytauE7}
[x_\sigma,y_\tau]=\epse{7}(\sigma,\tau)
\varphi_{\sigma,\tau}(x_\sigma,y_\tau)\,,
\end{equation}
for any $\sigma,\tau\in \calSe{7}$, $x_\sigma\in V(\sigma)$ and
$y_\tau\in V(\tau)$; where
$\epse{7}:\calSe{7}\times\calSe{7}\rightarrow F$ is given by Table
\ref{ta:e7}.

\begin{table}[h!]
{\footnotesize
$$
\vbox{\offinterlineskip \halign{\hfil$#$\hfil\vreglon
 &\hfil$#$\hfil\vregla
 &\hfil$#$\hfil\vregla
 &\hfil$#$\hfil\vregla
 &\hfil$#$\hfil&\hskip .2pt\vregla#
 &&\hfil$#$\hfil\vregla\cr
 \omit\enspace\hfil\raise2pt\hbox{$\epse{7}$}%
  \enspace\hfil\vrule width1pt depth 4pt
   &\quad\emptyset\quad&\{1,2,5,6\}&\{2,3,6,7\}&\{1,3,5,7\}&\omit\vregleta
     &$\{1,2,3,4\}$&\{3,4,5,6\}&\{1,4,6,7\}&\{2,4,5,7\}\cr
 \noalign{\hreglon}
 \emptyset&1&1&1&1&\omit\vregleta&1&1&1&1\cr
 \noalign{\hregla}
 \{1,2,5,6\}&1&-1&1&-1&\omit\vregleta&1&-1&-1&1\cr\noalign{\hregla}
 \{2,3,6,7\}&1&-1&-1&1&\omit\vregleta&1&1&-1&-1\cr\noalign{\hregla}
 \{1,3,5,7\}&1&1&-1&-1&\omit\vregleta&1&-1&1&-1\cr
 \multispan6{\hregletafill}&\multispan4{\hreglafill}\cr
 \{1,2,3,4\}&1&-1&-1&-1&&-1&1&1&1\cr\noalign{\hregla}
 \{3,4,5,6\}&1&1&-1&1&&-1&-1&-1&1\cr\noalign{\hregla}
 \{1,4,6,7\}&1&1&1&-1&&-1&1&-1&-1\cr\noalign{\hregla}
 \{2,4,5,7\}&1&-1&1&1&&-1&-1&1&-1\cr
 \noalign{\hregla}
}}
$$}\caption{$\epse{7}$}\label{ta:e7}
\end{table}

All the entries that appear in this table are $\pm 1$, and these
signs coincide with the signs that appear in the multiplication
table of the real algebra of octonions $\bO$: The usual basis of
this algebra is $\{1,i,j,k(=ij),l,il,jl,kl\}$, with multiplication
table given in Table \ref{ta:octonions}.

\begin{table}[h!]
{\small
$$
\vbox{\offinterlineskip \halign{\hfil$#$\enspace\hfil\vreglon
 &\hfil\enspace$#$\enspace\hfil\vregla
 &\hfil\enspace$#$\enspace\hfil\vregla
 &\hfil\enspace$#$\enspace\hfil\vregla
 &\hfil\enspace$#$\enspace\hfil&\hskip .2pt\vregla#
 &&\hfil\enspace$#$\enspace\hfil\vregla\cr
 \omit\hfil\vrule width1pt depth 4pt
   &1&i&j&k&\omit\vregleta
     &l&il&jl&kl\cr
 \noalign{\hreglon}
 1&1&i&j&k&\omit\vregleta&l&il&jl&kl\cr
 \noalign{\hregla}
 i&i&-1&k&-j&\omit\vregleta&il&-l&-kl&jl\cr\noalign{\hregla}
 j&j&-k&-1&i&\omit\vregleta&jl&kl&-l&-il\cr\noalign{\hregla}
 k&k&j&-i&-1&\omit\vregleta&kl&-jl&il&-l\cr
 \multispan6{\hregletafill}&\multispan4{\hreglafill}\cr
 l&l&-il&-jl&-kl&&-1&i&j&k\cr\noalign{\hregla}
 il&il&l&-kl&jl&&-i&-1&-k&k\cr\noalign{\hregla}
 jl&jl&kl&l&-il&&-j&k&-1&-i\cr\noalign{\hregla}
 kl&kl&-jl&il&l&&-k&-j&i&-1\cr
 \noalign{\hregla}
}}
$$}\caption{The Octonions}\label{ta:octonions}
\end{table}

Actually, $\fre_7=\frg(S_8,S_4)$ is $\bZ_2\times\bZ_2$-graded, and
$S_8$ is naturally $\bZ_2$-graded (with even part $S_4$
\eqref{eq:bullet}). Thus $\fre_7$ is
$\bZ_2\times\bZ_2\times\bZ_2$-graded. In fact, the decomposition
\eqref{eq:e7} is a `fake' $\bZ_2^7$-grading, as the $\bZ_2$-linear
map
\begin{equation}\label{eq:chi3}
\begin{split}
\chi_3:\bZ_2^3\quad &\longrightarrow\quad \bZ_2^7\\
(x_1,x_2,x_3)&\mapsto
(x_1+x_3,x_1+x_2+x_3,x_2+x_3,x_3,x_1,x_1+x_2,x_2)
\end{split}
\end{equation}
provides a bijection between $\bZ_2^3$ and $\calSe{7}\subseteq
\bZ_2^7$. Through $\chi_3$, we get the $\bZ_2^3$-grading of
$\fre_7$.

Following \cite{HelenaMajid}, $\bO$ is a \emph{twisted group
algebra} $\bO=\bR_\phi[\bZ_2^3]$. This is the algebra defined on
the group algebra $\bR[\bZ_2^3]$ (with basis $\{ e^x: x\in
\bZ_2^3\}$ and multiplication determined by $e^xe^y=e^{x+y}$), but
with a new multiplication of the basic elements: $e^x\cdot
e^y=\phi(x,y)e^{x+y}$ for any $x,y\in\bZ_2^3$, where
 $\phi(x,y)=(-1)^{f(x,y)}$ and
$f:\bZ_2^3\times\bZ_2^3\rightarrow \bZ_2$ is given by
\begin{equation}\label{eq:f}
f(x,y)=\bigl(\sum_{i\geq j}x_iy_j\bigr)
+x_1y_2y_3+x_2y_1y_3+x_3y_1y_2\,,
\end{equation}
for $x=(x_1,x_2,x_3),y=(y_1,y_2,y_3)$ in $\bZ_2^3$. (There is a
discrepancy here with respect to \cite{HelenaMajid} due to the
fact that these authors work with the basis
$\{1,i,j,k,l,li,lj,lk\}$.) One can check easily that identifying
$i$ with $(1,0,0)$, $j$ with $(0,1,0)$ and $k$ with $(0,0,1)$ one
recovers Table \ref{ta:octonions}.

Note that through $\chi_3$, this identifies $i$ with
$\{1,2,5,6\}$, $j$ with $\{2,3,6,7\}$ and $k$ with $\{1,2,3,4\}$.
Moreover $\calSe{7}=\chi_3\bigl(\bZ_2^3\bigr)$ and
\begin{equation}\label{eq:epsilonf}
\epse{7}(\sigma,\tau)=
(-1)^{f\bigl(\chi_3^{-1}(\sigma),\chi_3^{-1}(\tau)\bigr)}
\end{equation}
for any $\sigma,\tau\in \calSe{7}$, thus providing a closed
formula for $\epse{7}$.

\bigskip

Also, the split Lie algebra of type $D_6$ is realized as
$\frd_6=\frg(S_4,S_4)$, a subalgebra of  $\fre_7=\frg(S_8,S_4)$.
Thus, by restriction, one gets
\begin{equation}\label{eq:d6}
\frd_6=\oplus_{\sigma\in\calS_{\frd_6}}V(\sigma)\,,
\end{equation}
where
\[
\calS_{\frd_6}=\Bigl\{\emptyset,\{1,2,5,6\},\{2,3,6,7\},\{1,3,5,7\}\Bigr\}\subseteq
2^{\{1,2,3,5,6,7\}}.
\]
and the multiplication is given by
\begin{equation}\label{eq:xsigmaytauD6}
[x_\sigma,y_\tau]=\epsilon_{\frd_6}(\sigma,\tau)
\chi_{\sigma,\tau}(x_\sigma,y_\tau)\,,
\end{equation}
for any $\sigma,\tau\in \calS_{\frd_6}$, $x_\sigma\in V(\sigma)$
and $y_\tau\in V(\tau)$; where
$\epsilon_{\frd_6}:\calS_{\frd_6}\times\calS_{\frd_6}\rightarrow
F$ is given by the left top corner of Table \ref{ta:e7}.
Therefore, we may consider the $\bZ_2$-linear map:
\begin{equation}\label{eq:chi2}
\begin{split}
\chi_2:\bZ_2^2\quad &\longrightarrow\quad \bZ_2^6\,\bigl(\cong 2^{\{1,2,3,5,6,7\}}\bigr)\\
(x_1,x_2)&\mapsto (x_1,x_1+x_2,x_2,x_1,x_1+x_2,x_2)
\end{split}
\end{equation}
and the map $f':\bZ_2^2\times\bZ_2^2\rightarrow \bZ_2$ given by
\[
f'\bigl((x_1,x_2),(y_1,y_2)\bigr)=x_1y_1+x_2(y_1+y_2) \,.
\]
Then, $\calS_{\frd_6}=\chi_2\bigl(\bZ_2^2\bigr)$ and
\begin{equation}\label{eq:epsilonf'}
\epsilon_{\frd_6}(\sigma,\tau)=
(-1)^{f'\bigl(\chi_2^{-1}(\sigma),\chi_2^{-1}(\tau)\bigr)}
\end{equation}
for any $\sigma,\tau\in \calS_{\frd_6}$, thus providing a closed
formula too for $\epsilon_{\frd_6}$.

\medskip

\subsection{$\mathbf{E_8}$} \quad\null

All the arguments for $E_7$ can be easily extended to $E_8$. The
split simple Lie algebra $E_8$ is realized as
$\fre_8=\frg(S_8,S_8)$. Here there are $4$ copies of $V$ involved
in the first $S_8$ and another $4$ copies for the second $S_8$.
The indices $1,2,3,4$ will be used for the first $S_8$ and the
indices $5,6,7,8$ for the second $S_8$. Therefore, with the same
arguments as above,
\begin{equation}\label{eq:e8}
\fre_8=\oplus_{\sigma\in\calSe{8}}V(\sigma)\,,
\end{equation}
where
\[
\begin{split}
\calSe{8}&=\Bigl\{\emptyset,\{1,2,5,6\},\{2,3,6,7\},\{1,3,5,7\},\\
  &\qquad \{1,2,3,4\},\{3,4,5,6\},\{1,4,6,7\},\{2,4,5,7\},\\
  &\qquad \{3,4,7,8\},\{1,4,5,8\},\{2,4,6,8\},\{5,6,7,8\},\\
  &\qquad \{1,2,7,8\},\{2,3,5,8\},\{1,3,6,8\}\Bigr\}\\
  &\subseteq 2^{\{1,2,3,4,5,6,7,8\}}.
\end{split}
\]
and the multiplication is given by
\begin{equation}\label{eq:xsigmaytauE8}
[x_\sigma,y_\tau]=\epse{8}(\sigma,\tau)
\varphi_{\sigma,\tau}(x_\sigma,y_\tau)\,,
\end{equation}
for any $\sigma,\tau\in \calSe{8}$, $x_\sigma\in V(\sigma)$ and
$y_\tau\in V(\tau)$; where
$\epse{8}:\calSe{8}\times\calSe{8}\rightarrow F$ is given by Table
\ref{ta:e8}.

\begin{table}[h!]
{\footnotesize
$$
\vbox{\offinterlineskip \halign{\hfil$#$\hfil\vreglon
 &&\hfil$#$\hfil\vregla\cr
 \omit\enspace\hfil\raise2pt\hbox{$\epse{8}$}%
 \enspace\hfil\vrule width1pt
   &\enspace\emptyset\enspace
   &\begin{Bmatrix}1\\2\\5\\6\end{Bmatrix}
   &\begin{Bmatrix}2\\3\\6\\7\end{Bmatrix}
   &\begin{Bmatrix}1\\3\\5\\7\end{Bmatrix}
   &\begin{Bmatrix}1\\2\\3\\4\end{Bmatrix}
   &\begin{Bmatrix}3\\4\\5\\6\end{Bmatrix}
   &\begin{Bmatrix}1\\4\\6\\7\end{Bmatrix}
   &\begin{Bmatrix}2\\4\\5\\7\end{Bmatrix}
   &\quad
   &\begin{Bmatrix}3\\4\\7\\8\end{Bmatrix}
   &\begin{Bmatrix}1\\4\\5\\8\end{Bmatrix}
   &\begin{Bmatrix}2\\4\\6\\8\end{Bmatrix}
   &\begin{Bmatrix}5\\6\\7\\8\end{Bmatrix}
   &\begin{Bmatrix}1\\2\\7\\8\end{Bmatrix}
   &\begin{Bmatrix}2\\3\\5\\8\end{Bmatrix}
   &\begin{Bmatrix}1\\3\\6\\8\end{Bmatrix}\cr
 \noalign{\hreglon}
 \emptyset&1&1&1&1&1&1&1&1
   &&1&1&1&1&1&1&1\cr
 \noalign{\hregla}
 \{1,2,5,6\}&1&-1&1&-1&1&-1&-1&1&
 &&1&-1&1&-1&-1&1\cr\noalign{\hregla}
 \{2,3,6,7\}&1&-1&-1&1&1&1&-1&-1&
 &-1&&1&1&1&-1&-1\cr\noalign{\hregla}
 \{1,3,5,7\}&1&1&-1&-1&1&-1&1&-1&
 &1&-1&&1&-1&1&-1\cr\noalign{\hregla}
 \{1,2,3,4\}&1&-1&-1&-1&-1&1&1&1&
 &1&1&1&&-1&-1&-1\cr\noalign{\hregla}
 \{3,4,5,6\}&1&1&-1&1&-1&-1&-1&1&
 &-1&1&-1&1&&1&-1\cr\noalign{\hregla}
 \{1,4,6,7\}&1&1&1&-1&-1&1&-1&-1&
 &-1&-1&1&1&-1&&1\cr\noalign{\hregla}
 \{2,4,5,7\}&1&-1&1&1&-1&-1&1&-1&
 &1&-1&-1&1&1&-1&\cr\noalign{\hregla}
 \omit\hfil\vrule height 10pt width1pt &&&&&&&&&
 &&&&&&&\cr\noalign{\hregla}
 \{3,4,7,8\}&1&&1&-1&-1&1&1&-1&
 &-1&1&-1&-1&1&1&-1\cr\noalign{\hregla}
 \{1,4,5,8\}&1&-1&&1&-1&-1&1&1&
 &-1&-1&1&-1&-1&1&1\cr\noalign{\hregla}
 \{2,4,6,8\}&1&1&-1&&-1&1&-1&1&
 &1&-1&-1&-1&1&-1&1\cr\noalign{\hregla}
 \{5,6,7,8\}&1&-1&-1&-1&&-1&-1&-1&
 &1&1&1&-1&1&1&1\cr\noalign{\hregla}
 \{1,2,7,8\}&1&1&-1&1&1&&1&-1&
 &-1&1&-1&-1&-1&-1&1\cr\noalign{\hregla}
 \{2,3,5,8\}&1&1&1&-1&1&-1&&1&
 &-1&-1&1&-1&1&-1&-1\cr\noalign{\hregla}
 \{1,3,6,8\}&1&-1&1&1&1&1&-1&&
 &1&-1&-1&-1&-1&1&-1\cr\noalign{\hregla} }}
$$}\caption{$\epse{8}$}\label{ta:e8}
\end{table}

The $\pm$ signs that appear in Table \ref{ta:e8} are the same as
those that appear in the corresponding entry of the multiplication
table of $\bO\oplus\bO=\bO\otimes_{\bR}\bR[\varepsilon]$ (with
$\bR[\varepsilon]=\bR 1\oplus\bR\varepsilon$ and
$\varepsilon^2=1$) in the basis
\[
\begin{split}
\{1\otimes 1&,i\otimes 1,j\otimes 1,k\otimes 1,l\otimes
1,(il)\otimes 1,(jl)\otimes 1,(kl)\otimes 1,\\
& 1\otimes\varepsilon,
i\otimes\varepsilon,j\otimes\varepsilon,k\otimes\varepsilon,
-l\otimes\varepsilon,-(il)\otimes\varepsilon,-(jl)\otimes\varepsilon,
-(kl)\otimes\varepsilon\}
\end{split}
\]
(Note the minus signs of the last four elements.) Working with
this basis, $\bO\otimes_\bR \bR[\varepsilon]$ is shown to be
isomorphic to the twisted group algebra $\bR_\Phi[\bZ_2^4]$, with
$\Phi(x,y)=(-1)^{g(x,y)}$ for $x,y\in\bZ_2^4$, where
\[
\begin{split}
g\bigl((x_1,x_2,x_3,&x_4),(y_1,y_2,y_3,y_4)\bigr)\\&=
 f\bigl((x_1,x_2,x_3),(y_1,y_2,y_3)\bigr)+x_3y_4+x_4y_3\\
 &=\bigl(\sum_{1\leq j\leq i\leq 3}x_iy_j\bigr)
 +x_3y_4+x_4y_3+x_1y_2y_3+x_2y_1y_3+x_3y_1y_2\,,
\end{split}
\]
with $f$ as in \eqref{eq:f}. The isomorphism carries $i\otimes 1$
to $(1,0,0,0)$, $j\otimes 1$ to $(0,1,0,0)$, $l\otimes 1$ to
$(0,0,1,0)$ and $1\otimes\varepsilon$ to $(0,0,0,1)$; so that, for
instance, $l\otimes\varepsilon$ goes to
\[
(-1)^{g\bigl((0,0,1,0),(0,0,0,1)\bigr)}(0,0,1,1)=-(0,0,1,1).
\]

As for $\fre_7$, consider the $\bZ_2$-linear map
\begin{equation}\label{eq:chi4}
\begin{split}
\chi_4:\bZ_2^4\ &\longrightarrow\ \bZ_2^8\\
(x_1,x_2,x_3,x_4)&\mapsto
(x_1+x_3+x_4,x_1+x_2+x_3+x_4,x_2+x_3+x_4,\\
&\phantom{\longrightarrow}\quad
x_3+x_4,x_1+x_4,x_1+x_2+x_4,x_2+x_4,x_4)
\end{split}
\end{equation}
which satisfies that $\calSe{8}=\chi_4\bigl(\bZ_2^4\setminus
\{(0,0,0,1)\}\bigr)$ and
\begin{equation}\label{eq:epsilong}
\epse{8}(\sigma,\tau)=
(-1)^{g\bigl(\chi_4^{-1}(\sigma),\chi_4^{-1}(\tau)\bigr)}
\end{equation}
for any $\sigma,\tau\in \calSe{8}$, which gives a closed formula
for $\epse{8}$.

\bigskip

\begin{remark} The models thus constructed for the split $E_7$ and
$E_8$ are strongly related to a very interesting combinatorial
construction previously given by A.~Grishkov in \cite{Grishkov}.
To see this, take a basis $\{v,w\}$ of $V$ with $\langle v\vert
w\rangle =1$ and denote by $\{v_i,w_i\}$ the copy of this basis in
$V_i$. Take, for instance, $\sigma=\{1,2,7,8\}$ and
$\tau=\{1,7,8\}\subseteq \sigma$; then the element $(\sigma,\tau)$
in Grishkov's construction correspond to the element $v_1\otimes
w_2\otimes v_7\otimes v_8$ in our $V(\sigma)\subseteq \fre_8$. It
is hoped that this will make Grishkov's construction to appear
more natural. In \cite{Grishkov} one has to choose an
identification of the Moufang loop $M_7$ consisting of the
elements $\{\pm1,\pm i,\pm j\pm k,\pm l,\pm il,\pm jl,\pm kl\}$ in
$\bO$ with $\pm\calSe{7}$ and another of the Moufang loop
$M_8=M_7\times\bZ_2$ with
$\pm\Bigl(\calSe{8}\cup\bigl\{\{1,2,3,4,5,6,7,8\}\bigr\}\Bigr)$,
but no hint is given as  how to perform this subtle point. Here
such identifications have been explicitly given.
\end{remark}

\begin{remark}\label{re:Fano}
The referee has kindly suggested nice interpretations of
$\calSe{7}\setminus\{\emptyset\}$ as the set of quadrilaterals of
the Fano projective plane $\bZ_2\bP^2$, and of
$\calSe{8}\setminus\{\emptyset\}$ as the set of affine planes in
the affine space $\bZ_2^3$, once the set $\{1,2,3,4,5,6,7\}$ is
suitably identified with the set of points of the Fano plane, and
$\{1,2,3,4,5,6,7,8\}$ with $\bZ_2^3$.

For instance, identify the set $\{1,2,3,4,5,6,7\}$ with the set of
points of the Fano plane (i.e., the nonzero elements in $\bZ_2^3$)
as follows:
\begin{equation}\label{eq:unosiete}
\begin{aligned}
1&\leftrightarrow (1,0,1)& 2&\leftrightarrow (1,1,1)&
3&\leftrightarrow (0,1,1)& 4&\leftrightarrow (0,0,1)\\
5&\leftrightarrow (1,0,0)& 6&\leftrightarrow (1,1,0)&
7&\leftrightarrow (0,1,0)& &
\end{aligned}
\end{equation}
Then the elements of $\calSe{7}\setminus\{\emptyset\}$ are
precisely the sets of vertices of the quadrilaterals (the
complements of the lines) in the Fano plane. By projective
duality, these quadrilateral (or better, their complementary
lines) are in bijection with the set of points of the Fano plane,
and this leads to the bijection of $\calSe{7}\setminus
\{\emptyset\}$ and $\bZ_2^3\setminus\{(0,0,0)\}$ given by $\chi_3$
in \eqref{eq:chi3}.

Also, the set $\{1,2,3,4,5,6,7,8\}$ may be identified with
$\bZ_2^3$ by means of \eqref{eq:unosiete}, together with
$8\leftrightarrow (0,0,0)$. If $l$ is a line in $\bZ_2\bP^2\simeq
\bZ_2^3\setminus\{\emptyset\}$, then $l\cup\{(0,0,0)\}$ is an
affine plane in $\bZ_2^3$ through the origin. On the other hand,
if $q$ is a quadrilateral in $\bZ_2\bP^2$, then $q$ is an affine
plane in $\bZ_2^3$ which does not contain the origin. The elements
in $\calSe{8}\setminus\{\emptyset\}$ are precisely the affine
planes in $\bZ_2^3$, and with our identification
\eqref{eq:unosiete}, for any $(a_1,a_2,a_3,a_4)\in \bZ_2^4$,
$\chi_4(a_1,a_2,a_3,a_4)$, considered as a subset of
$\{1,2,3,4,5,6,7,8\}$, and hence of $\bZ_2^3$, is precisely the
affine plane with equation $a_1x+a_2y+a_3z+a_4=1$.
\end{remark}

\medskip

\subsection{$\mathbf{E_6}$} \quad\null

For $E_6$ things are a bit more involved, partly because it does
not contain a subalgebra isomorphic to $\frsp(V)^6$. Let
$K=F\times F$ be the split two dimensional unital composition
algebra and let $S_2$ be the associated para-Hurwitz algebra. Thus
$S_2$ has a basis $\{e^+,e^-\}$ with multiplication given by
\[
e^\pm\bullet e^\pm =e^\mp ,\qquad e^\pm\bullet e^\mp=0\,,
\]
and with norm given by $q(e^\pm)=0$ and $q(e^+,e^-)=1$. The
orthogonal Lie algebra $\fro(S_2,q)$ is spanned by
$\phi=\sigma_{e^-,e^+}$, which satisfies $\phi(e^\pm)=\pm e^\pm$,
and the triality Lie algebra is
\[
\tri(S_2,\bullet,q)=\bigl\{ (\alpha\phi,\beta\phi,\gamma\phi) :
 \alpha,\beta,\gamma\in F,\ \alpha+\beta+\gamma=0\bigr\}\,.
\]
Besides,
\[
\begin{split}
t_{e^-,e^+}&=
 \Bigl(\sigma_{e^-,e^+},\frac{1}{2}q(e^-,e^+)id -r_{e^-}l_{e^+},
    \frac{1}{2}q(e^-,e^+)id -l_{e^-}r_{e^+}\Bigr)\\
 &=\Bigl(\phi,-\frac{1}{2}\phi,-\frac{1}{2}\phi\Bigr)\,.
\end{split}
\]
Let $\frt_2=\{(\alpha_0,\alpha_1,\alpha_2)\in
k^3:\alpha_0+\alpha_1+\alpha_2=0\}$, which is a two-dimensional
abelian Lie algebra. For any $\sigma\in \calSf\setminus
\Bigl\{\emptyset,\{1,2,3,4\}\Bigr\}$, consider the element
$a_\sigma\in \frt_2$ given by $a_\sigma=a_{\bar\sigma}$, where
$\bar\sigma=\{1,2,3,4\}\setminus\sigma$, and
\[
a_{\{ 1,2\}}=(1,-\tfrac{1}{2},-\tfrac{1}{2}), \quad
 a_{\{ 2,3\}}=(-\tfrac{1}{2},1,-\tfrac{1}{2}),\quad
 a_{\{ 1,3\}}=(-\tfrac{1}{2},-\tfrac{1}{2},-1).
\]
Consider too the $\sigma$-action of $\frt_2$ on $E=S_2=Fe^++Fe^-$
defined by
\[
(\alpha_0,\alpha_1,\alpha_2).e^\pm =\alpha_i e^\pm,
\]
where $i=0$ for $\sigma=\{1,2\}$ or $\{3,4\}$, $i=1$ for
$\sigma=\{2,3\}$ or $\{1,4\}$, and $i=2$ for $\sigma=\{1,3\}$ or
$\{2,4\}$.


Then, the split Lie algebra $E_6$ is realized as $\frg(S_8,S_2)$
which, according to \eqref{eq:gss'}, is described as:
\begin{equation}\label{eq:e6}
\fre_6=\oplus_{\sigma\in\calSf}\tilde V(\sigma)\,,
\end{equation}
where
\[
\begin{cases}
\tilde V(\emptyset)=\frsp(V)^4\oplus\frt_2\\
\tilde V\bigl(\{1,2,3,4\}\bigr)
  =V\bigl(\{1,2,3,4\})\,\bigl(=V_1\otimes
  V_2\otimes V_3\otimes V_4\bigr)\\
\tilde V\bigl(\sigma)=V(\sigma)\otimes E\quad\text{for any
$\sigma\in\calSf$, $\sigma\ne\emptyset,\{1,2,3,4\}$.}
\end{cases}
\]
The vector space $\tilde V(\sigma)$ is a module for
$\frsp(V)^4\oplus\frt_2$ by means of the natural action of
$\frsp(V)^4$ on $V(\sigma)$ and the $\sigma$-action of $\frt_2$ on
$E$ if $\sigma\in\calSf\setminus
\Bigl\{\emptyset,\{1,2,3,4\}\Bigr\}$.

Finally, take the bilinear maps $\varphi_{\sigma,\tau}$ defined in
\eqref{eq:phisigmatau} and define the new bilinear maps
\[
\tilde\varphi_{\sigma,\tau}: \tilde V(\sigma)\times \tilde
V(\tau)\rightarrow \tilde V(\sigma+\tau)\,,
\]
for $\sigma,\tau\in\calSf$, as follows:

\begin{itemize}

\item $\tilde\varphi_{\emptyset,\emptyset}$ is the Lie bracket in
the direct sum $\frsp(V)^4\oplus\frt_2$.
\item
 For any $\emptyset\ne\sigma\in\calSf$,
 $\tilde\varphi_{\emptyset,\sigma}=-\tilde\varphi_{\sigma,\emptyset}$
 is given by the action of $\frsp(V)^4\oplus\frt_2$ on $\tilde
 V(\sigma)$. ($\frsp(V)^4$ acts on $V(\sigma)$ and, if present, $\frt_2$ on
 $E$ by means of the $\sigma$-action.)
\smallskip
\item
 $\tilde\varphi_{\{1,2,3,4\},\{1,2,3,4\}}
 =\varphi_{\{1,2,3,4\},\{1,2,3,4\}}$.
\smallskip
\item
 For any $\sigma\in\calSf\setminus
  \Bigl\{\emptyset,\{1,2,3,4\}\Bigr\}$,
\[
\begin{split}
\tilde\varphi_{\{1,2,3,4\},\sigma}
  (x_{\{1,2,3,4\}},y_\sigma\otimes e)&=
  \varphi_{\{1,2,3,4\},\sigma}
  (x_{\{1,2,3,4\}},y_\sigma)\otimes e\\
  &=\tilde\varphi_{\sigma,\{1,2,3,4\}}
  (y_\sigma\otimes e,x_{\{1,2,3,4\}})
\end{split}
\]
for any $x_{\{1,2,3,4\}}\in V(\{1,2,3,4\})=\tilde V(\{1,2,3,4\})$,
$y_\sigma\in V(\sigma)$ and $e\in E$.
\smallskip
\item For any $\sigma\in\calSf\setminus
  \Bigl\{\emptyset,\{1,2,3,4\}\Bigr\}$, so $\sigma=\{i,j\}$ for some
  $i,j\in\{1,2,3,4\}$, and any $v_i,w_i\in V_i$, $v_j,w_j\in V_j$,
  and $\nu,\nu'\in\{+,-\}$:
\[
\begin{split}
\tilde\varphi_{\sigma,\sigma}&
  \bigl(v_i\otimes v_j\otimes e^\nu,w_i\otimes w_j\otimes e^{\nu'}\bigr)\\
  &=\delta_{\nu,-\nu'}
   \Bigl(\varphi_{\sigma,\sigma}(v_i\otimes v_j,w_i\otimes w_j)
     -\nu \langle v_i\vert w_i\rangle\langle v_j\vert w_j\rangle
     a_\sigma\Bigr)\, ,
\end{split}
\]
which belongs to $\frsp(V)^4\oplus\frt_2=\tilde V(\emptyset)$,
where $\delta$ is the Kronecker delta; while, with
$\bar\sigma=\{1,2,3,4\}\setminus \sigma$ as before,
\[
\null\qquad \tilde\varphi_{\sigma,\bar\sigma}
 \bigl(x_\sigma\otimes e^\nu,y_{\bar\sigma}\otimes e^{\nu'}\bigr)
 =\delta_{\nu,-\nu'}\varphi_{\sigma,\bar\sigma}(x_\sigma,y_{\bar\sigma})
 \in \tilde V(\{1,2,3,4\}),
\]
for any $x_\sigma\in V(\sigma)$, $y_{\bar\sigma}\in
V(\bar\sigma)$.
\smallskip
\item Finally, for any $\sigma,\tau\in\calSf\setminus
  \Bigl\{\emptyset,\{1,2,3,4\}\Bigr\}$ with
  $\tau\ne\sigma,\bar\sigma$,
\[
\tilde\varphi_{\sigma,\tau}
 \bigl(x_\sigma\otimes e^\nu,y_{\tau}\otimes e^{\nu'}\bigr)
 =\delta_{\nu,\nu'}\varphi_{\sigma,\tau}(x_\sigma,y_\tau)\otimes
 e^{-\nu}
\]
for any $x_\sigma\in V(\sigma)$, $y_\tau\in V(\tau)$ and
$\nu,\nu'\in \{+,-\}$.
\end{itemize}

\bigskip

Again, the multiplication in $\fre_6=\frg(S_8,S_2)$ is given by
\begin{equation}\label{eq:xsigmaytauE6}
[x_\sigma,y_\tau]=\epse{6}(\sigma,\tau)\tilde\varphi(x_\sigma,y_\tau)
\,,
\end{equation}
for any $\sigma,\tau\in \calSf$, $x_\sigma\in \tilde V(\sigma)$,
$y_\tau\in \tilde V(\tau)$, where
\[
\epse{6}(\sigma,\tau)=
\begin{cases}
\frac{1}{2}\epsf (\sigma,\tau)&\text{if
$\sigma,\tau\in\calSf\setminus
  \Bigl\{\emptyset,\{1,2,3,4\}\Bigr\}$,
  $\tau=\sigma$ or $\tau=\bar\sigma$,}\\
\epsf(\sigma,\tau)&\text{otherwise.}
\end{cases}
\]
\medskip

\begin{remark}\label{re:char3}
Note that in characteristic $3$,
$a_{\{1,2\}}=a_{\{2,3\}}=a_{\{1,3\}}=(1,1,1)$, and hence
${\hat\fre}_6=[\fre_6,\fre_6]=\oplus_{\sigma\in\calSf}\hat
V(\sigma)$, with $\hat V(\sigma)=\tilde V(\sigma)$ for any
$\sigma\ne\emptyset$, but with $\hat
V(\emptyset)=\frsp(V)^4\oplus{\hat\frt}_2$, where
${\hat\frt}_2=Fa_{\{1,2\}}$. Actually, the Lie algebra obtained by
taking the $\bZ$-span of a Chevalley basis of the complex simple
Lie algebra of type $E_6$, and then tensoring with a field of
characteristic three, is no longer simple, but has a one
dimensional center (see, for instance, \cite{D}, or \cite[\S
3]{VK}). Modulo this center, one has a simple Lie algebra of
dimension $77$, which is isomorphic to $\hat\fre_6$, and $\fre_6$
is isomorphic to the Lie algebra of derivations
$\Der(\hat\fre_6)$.
\end{remark}

\bigskip

Also, the split Lie algebra $A_5$ is realized as
$\fra_5=\frg(S_4,S_2)$, a subalgebra of $\fre_6=\frg(S_8,S_2)$.
Thus, by restriction, one gets
\begin{equation}\label{eq:a5}
\fra_5=\oplus_{\sigma\in\calS_{\frc_3}}\hat V(\sigma)\,,
\end{equation}
where $\hat V(\emptyset)=\frsp(V)^3\oplus\frt_2$ and $\hat
V(\sigma)=\tilde V(\sigma)$ for any $\emptyset\ne\sigma\in
\calS_{\frc_3}$, and the multiplication is given by
\begin{equation}\label{eq:xsigmaytauA5}
[x_\sigma,y_\tau]=\epsilon_{\fra_5}(\sigma,\tau)\tilde\varphi(x_\sigma,y_\tau)
\,,
\end{equation}
for any $\sigma,\tau\in \calS_{\frc_3}$, $x_\sigma\in \hat
V(\sigma)$, $y_\tau\in \hat V(\tau)$, with
\[
\epsilon_{\fra_5}(\sigma,\tau)=
\begin{cases}
\frac{1}{2}\epsilon_{\frc_3}(\sigma,\tau)&\text{if
$\sigma=\tau\ne\emptyset$,}\\
\epsilon_{\frc_3}(\sigma,\tau)&\text{otherwise.}
\end{cases}
\]

Also here, if the characteristic is $3$, $\fra_5$ is no longer
simple, since it contains the simple codimension one ideal $\hat
\fra_5=[\fra_5,\fra_5]$, which is isomorphic to the projective
special linear algebra $\frpsl_6(F)$. It can be checked that
$\fra_5$ is isomorphic to the projective general linear algebra
$\frpgl_6(F)$.

\bigskip

\begin{remark}\label{re:rootsystem}
In concluding this section, let us remark that the  models of the
exceptional simple Lie algebras obtained provide, in particular,
very easy descriptions of the exceptional root systems, different
from the description in \cite{Bou2}. Thus, for instance, for
$E_8$, take a symplectic basis $\{ v_i,w_i\}$ of $V_i$
($i=1,\ldots,8$) (that is, $\langle v_i\vert w_i\rangle =1$). Then
the vector space $\frh_8=\espan{ \gamma_{v_i,w_i}: i=1,\ldots,8}\,
\bigl(\subseteq V(\emptyset)=\oplus_{i=1}^8\frsp(V_i)\bigr)$ is a
Cartan subalgebra of $\fre_8$. Since $\gamma_{v_i,w_i}(v_i)=-v_i$,
$\gamma_{v_i,w_i}(w_i)=w_i$, and hence also
$[\gamma_{v_i,w_i},\gamma_{v_i,v_i}]=-2\gamma_{v_i,v_i}$ and
$[\gamma_{v_i,w_i},\gamma_{w_i,w_i}]=2\gamma_{w_i,w_i}$, let
$\varepsilon_i:\frh_8\rightarrow F$ be the linear map given by
$\varepsilon_i(\gamma_{v_j,w_j})=\delta_{ij}$ for any $1\leq
i,j\leq 8$. The description of $\fre_8$ in \ref{eq:e8} shows that
the corresponding root system is
\begin{multline*}
\Phi=\{\pm 2\varepsilon_i:i=1,\ldots,8\}\\
  \cup\{\pm\varepsilon_{i_1}
\pm\varepsilon_{i_2}\pm\varepsilon_{i_3}\pm\varepsilon_{i_4}:\sigma=\{i_1,i_2,i_3,i_4\}\in
\calS_{\fre_8}\setminus\{\emptyset\}\}.
\end{multline*}
Alternatively (see Remark \ref{re:Fano}), if the points of the
affine space $\bZ_2^3$ are numbered from $1$ to $8$, then
\begin{multline*}
\Phi=\{\pm 2\varepsilon_i:i=1,\ldots,8\}\\
  \cup\{\pm\varepsilon_{i_1}
\pm\varepsilon_{i_2}\pm\varepsilon_{i_3}\pm\varepsilon_{i_4}:\sigma=\{i_1,i_2,i_3,i_4\}\
 \text{is an affine plane in $\bZ_2^3$}\}.
\end{multline*}
 Using the lexicographic order with
$0<\varepsilon_1<\varepsilon_2<\dots<\varepsilon_8$, the simple
system of roots is $\Delta=\{\alpha_1,\ldots,\alpha_8\}$
(numbering as in \cite{Bou2}) with
$\alpha_1=-\varepsilon_1-\varepsilon_4-\varepsilon_6+\varepsilon_7$,
$\alpha_2=2\varepsilon_2$, $\alpha_3=2\varepsilon_1$,
$\alpha_4=-\varepsilon_1-\varepsilon_2-\varepsilon_3+\varepsilon_4$,
$\alpha_5=2\varepsilon_3$,
$\alpha_6=-\varepsilon_3-\varepsilon_4-\varepsilon_5+\varepsilon_6$,
$\alpha_7=2\varepsilon_5$ and
$\alpha_8=-\varepsilon_5-\varepsilon_6-\varepsilon_7+\varepsilon_8$.
(For $\fre_7$ it is enough to suppress here $\varepsilon_8$ and
$\alpha_8$ above or, alternatively, the points of the Fano plane
can be numbered from $1$ to $7$, and the root system of $\fre_7$
becomes $\{\pm 2\varepsilon_i:i=1,\ldots,7\}
  \cup\{\pm\varepsilon_{i_1}
\pm\varepsilon_{i_2}\pm\varepsilon_{i_3}\pm\varepsilon_{i_4}:\sigma=\{i_1,i_2,i_3,i_4\}\
 \text{is a quadrilateral in $\bZ_2\bP^2$}\}$.)
\end{remark}

\medskip
\section{Freudenthal triple systems}

Freudenthal \cite{FRII,FRVIII} obtained the exceptional simple Lie
algebras in terms of some triple systems, later called Freudenthal
triple systems. Actually, if $T$ is such a system, the direct sum
of two copies of $T$ becomes a Lie triple system, and hence the
odd part of a $\bZ_2$-graded Lie algebra. In this section, the
exceptional split simple Freudenthal triple systems will be
recovered from the models given in the previous section of the
exceptional simple Lie algebras. For our purposes, the approach
given by Yamaguti and Asano \cite{YamAs} of Freudenthal's
construction is more suitable.

Let us
first recall some definitions and results.

\medskip

Let $T$ be a vector space endowed with a nonzero alternating
bilinear form $(.\vert.):T\times T\rightarrow F$, and a triple
product $T\times T\times T\rightarrow T$: $(x,y,z)\mapsto [xyz]$.
Then $\bigl(T,[...],(.\vert.)\bigr)$ is said to be a
\emph{symplectic triple system} (see \cite{YamAs}) if it satisfies
the following identitities:
\begin{subequations}\label{eq:STS}
\begin{align}
&[xyz]=[yxz]\label{eq:STSa}\\
&[xyz]-[xzy]=(x\vert z)y-(x\vert y)z+2(y\vert z)x\label{eq:STSb}\\
&[xy[uvw]]=[[xyu]vw]+[u[xyv]w]+[uv[xyw]]\label{eq:STSc}\\
&([xyu]\vert v)+(u\vert [xyv])=0\label{eq:STSd}
\end{align}
\end{subequations}
for any elements $x,y,z,u,v,w\in T$.

Note that \eqref{eq:STSb} can be written as
\begin{equation}\label{eq:STSbb}
[xyz]-[xzy]=\psi_{x,y}(z)-\psi_{x,z}(y)
\end{equation}
with $\psi_{x,y}(z)=(x\vert z)y+(y\vert z)x$. (The maps
$\psi_{x,y}$ span the symplectic Lie algebra $\frsp(T)$.)

Also, \eqref{eq:STSc} is equivalent to $d_{x,y}=[xy.]$ being a
derivation of the triple system. Let $\inder(T)$ be the linear
span of $\{d_{x,y}:x,y\in T\}$, which is a Lie subalgebra of
$\End(T)$. For any $x,y,z,a,b\in T$, by \eqref{eq:STSb} and
\eqref{eq:STSc},
\begin{equation}\label{eq:indersp}
\begin{split}
0&=[xy([zab]-[zba])]-(z\vert b)[xya]+(z\vert a)[xyb]-2(a\vert
b)[xyz]\\
&=\bigl( [[xyz]ab]-[[xyz]ba]\bigr) +
\bigl([z[xya]b]-[zb[xya]]\bigr) +
\bigl([za[xyb]]-[z[xyb]a]\bigr)\\
&\hspace{1in} -(z\vert b)[xya]+ (z\vert a)[xya] - 2(a\vert b)[xyz]\\
&=\bigl(([xyz]\vert b)+(z\vert [xyb])\bigr) a -
 \bigl(([xyz]\vert a)+(z\vert [xya])\bigr) b\\
 &\hspace{1in} +
 2\bigl( ([xya]\vert b)+(a\vert [xyb])\bigr) z .
 \end{split}
\end{equation}
Hence, if $\dim T\geq 3$, \eqref{eq:STSd} (that is,
$\inder(T)\subseteq \frsp(T)$) follows from \eqref{eq:STSb} and
\eqref{eq:STSc}. Also, with $a=z$ \eqref{eq:indersp} gives
\[
3\bigl( [xya]\vert b)+(a\vert [xyb])\bigr)=0,
\]
so the same applies if the characteristic of $F$ is $\ne 3$. This
was already noted, over fields of characteristic $0$, in
\cite{YamAs}. However, this is no longer true if the
characteristic of $F$ is $3$ and $\dim T=2$. First, notice that in
this case, the right hand side of \eqref{eq:STSb} becomes
\[
(x\vert z)y+(y\vert x)z+(z\vert y)x
\]
which is $0$, since it is skew symmetric on its arguments and
$\dim T=2$. Hence \eqref{eq:STSa} and \eqref{eq:STSb} merely say
that $[...]$ is symmetric on its arguments.

As a counterexample in characteristic $3$, take $T=Fa+Fb$ with
$(a\vert b)=1$ and $[...]$ determined by $[aaa]=[bbb]=[abb]=0$ and
$[aab]=b$. One checks easily that \eqref{eq:STSc} is satisfied,
but $d_{a,a}\not\in \frsp(T)$, so \eqref{eq:STSd} is not
satisfied.

\bigskip

Symplectic triple systems are strongly related to a particular
kind of $\bZ_2$-graded Lie algebras:

\begin{theorem}\label{th:STS}
Let $\bigl(T,[...],(.\vert .)\bigr)$ be a symplectic triple system
and let $\bigl(V,\langle.\vert.\rangle\bigr)$ be a two dimensional
vector space endowed with a nonzero alternating bilinear form.
Define the $\bZ_2$-graded algebra $\frg=\frgo\oplus\frguno$ with
\[
\begin{cases}
\frgo=\inder(T)\oplus \frsp(V)&\text{(direct sum of ideals)}\\
\frguno=T\otimes V
\end{cases}
\]
and anticommutative multiplication given by:
\begin{itemize}
\item
$\frgo$ is a Lie subalgebra of $\frg$,
\item
$\frgo$ acts naturally on $\frguno$; that is
\[
[d,x\otimes v]=d(x)\otimes v,\qquad [s,x\otimes v]=x\otimes s(v),
\]
for any $d\in\inder(T)$, $s\in \frsp(V)$, $x\in T$, and $v\in V$.
\item
For any $x,y\in T$ and $u,v\in V$:
\begin{equation}\label{eq:oddbracket}
[x\otimes u,y\otimes v]=\langle u\vert v\rangle d_{x,y}+
  (x\vert y)\gamma_{u,v}
\end{equation}
where $\gamma_{u,v}=\langle u\vert .\rangle v+\langle v\vert
.\rangle u$.
\end{itemize}
Then $\frg$ is a Lie algebra.

Conversely, given a $\bZ_2$-graded Lie algebra
$\frg=\frgo\oplus\frguno$ with
\[
\begin{cases}
\frgo=\frs\oplus \frsp(V)&\text{(direct sum of ideals),}\\
\frguno=T\otimes V&\text{(as a module for $\frgo$),}
\end{cases}
\]
where $T$ is a module for $\frs$, by $\frgo$-invariance of the Lie
bracket, equation \eqref{eq:oddbracket} is satisfied for an
alternating bilinear form $(.\vert .):T\times T\rightarrow F$ and
a symmetric bilinear map $d_{.,.}:T\times T\rightarrow \frs$.
Then, if $(.\vert .)$ is not $0$ and a triple product on $T$ is
defined by means of $[xyz]=d_{x,y}(z)$, $\bigl(T,[...],(.\vert
.)\bigr)$ is a symplectic triple system.
\end{theorem}
\begin{proof}
This is a reformulation and extension of the results in
\cite[p.~256]{YamAs}, suitable for our purposes. First, since
$[xy.]\in \frsp(T)$ for any $x,y\in T$, a straightforward
computation shows that the Jacobi identity holds in the algebra
$\frg$.

Conversely, if $\frg=\bigl(\frs\oplus\frsp(V)\bigr)\oplus
(T\otimes V)$ is a Lie algebra as above, \eqref{eq:STSa} follows
by the symmetry of $d_{x,y}$. Also, the Jacobi identity with an
element of $\frs$ and two basic tensors in $T\otimes V$ shows that
$(.\vert .)$ and the map $(x,y)\mapsto d_{x,y}$ are
$\frs$-invariant and, in particular, one gets \eqref{eq:STSc}.
Finally, the Jacobi identity $J(x\otimes w,y\otimes v,z\otimes
v)=0$ with $x,y,z\in T$ and $\langle v\vert w\rangle =1$ gives
\eqref{eq:STSb}.
\end{proof}

\medskip

Note that the Lie algebra $\frg$ in Theorem \ref{th:STS} is graded
by the root system BC$_1$ with grading subalgebra of type C$_1$
(see \cite{BeSm}).

\medskip

If $\{ v,w\}$ is a symplectic basis of $V$ (notation as in the
previous Theorem), so $\langle v\vert w\rangle=1$, then with
$h=-\gamma_{v,w}$, $e=\tfrac{1}{2}\gamma_{v,v}$ and
$f=-\tfrac{1}{2}\gamma_{w,w}$, $\{ h,e,f\}$ is a canonical basis
of $\frsp(V)\cong \frsl_2$ ($[h,e]=2e$, $[h,f]=-2f$, $[e,f]=h$)
with $h(v)=v$, $h(w)=-w$. Hence the Lie algebra $\frg$ above is
$5$-graded
\[
\frg=\frg_{-2}\oplus\frg_{-1}\oplus\frg_0\oplus\frg_1\oplus\frg_2,
\]
where $\frg_i=\{x\in \frg: [h,x]=ix\}$ for any $i$. Here
$\frg_2=Fe$, $\frg_{-2}=Ff$, $\frg_1=T\otimes v$,
$\frg_{-1}=T\otimes w$ and $\frg_0=\frs\oplus Fh$.

\medskip

Now, let $T$ be a vector space endowed with a nonzero alternating
bilinear form $(.\vert .):T\times T\rightarrow F$, and a triple
product $T\times T\times T\rightarrow T$, $(x,y,z)\mapsto xyz$.
Then $\bigl(T,xyz,(.\vert .)\bigr)$ is said to be a
\emph{Freudenthal triple system} (see \cite{Meyberg,Ferrar,Brown})
if it satisfies:
\begin{subequations}\label{eq:FT}
\begin{align}
&\text{$xyz$ is symmetric in its arguments,}\label{eq:FTa}\\
&\text{$( x\vert yzt)$ is symmetric in its
arguments,}\label{eq:FTb}\\
&(xyy)xz+(yxx)yz+(xyy\vert z)x+(yxx\vert z)y +(x\vert
z)xyy+(y\vert z)yxx=0 \label{eq:FTc}
\end{align}
\end{subequations}
for any $x,y,z,t\in T$.

Note that the bilinear form $(.\vert .)$ is not assumed to be
nondegenerate.

\begin{theorem}\label{th:FTS}
Let $(.\vert .)$ be an alternating bilinear form on the vector
space $T$ and let $xyz$ and $[xyz]$ be two triple products on $T$
related by $xyz=[xyz]-\psi_{x,y}(z)$ for any $x,y,z\in T$. Then
$\bigl(T,xyz,(.\vert .)\bigr)$ is a Freudenthal triple system if
and only if $\bigl(T,[xyz],(.\vert ,)\bigr)$ is a symplectic
triple system.
\end{theorem}
\begin{proof}
First assume that $\bigl( T,[...],(.\vert .)\bigr)$ is a
symplectic triple system and define $xyz=[xyz]-\psi_{x,y}(z)$ for
any $x,y,z\in T$. Then $xyz$ is symmetric in $x$ and $y$ since so
are $[xyz]$ and $\psi_{x,y}$, while \eqref{eq:STSbb} implies that
$xyz$ is symmetric in $y$ and $z$, thus proving \eqref{eq:FTa}.
Also, $(x\vert yzt)$ is symmetric in $y$, $z$ and $t$ and, since
both $[yz.]$ and $\psi_{y,z}$ belong to $\frsp(T)$, it follows
that $(x\vert yzt)$ is symmetric too in $x$ and $t$, proving
\eqref{eq:FTb}. Finally, with $d_{x,y}=[xy.]$ for any $x,y\in T$,
\[
\begin{split}
d_{x,xyy}&=d_{x,d_{x,y}(y)}-d_{x,\psi_{x,y}(y)}
     =d_{x,d_{x,y}(y)}-(x\vert y)d_{x,y},\\
d_{y,yxx}&=d_{y,d_{y,x}(x)}-d_{y,\psi_{y,x}(x)}
     =d_{y,d_{x,y}(x)}+(x\vert y)d_{x,y},
\end{split}
\]
so by \eqref{eq:STSc}
\[
d_{x,xyy}+d_{y,yxx}=d_{x,d_{x,y}(y)}+d_{d_{x,y}(x),y}=[d_{x,y},d_{x,y}]=0.
\]
Thus,
\[
d_{x,xyy}+d_{y,yxx}=0
\]
for any $x,y\in T$, and this is equivalent to \eqref{eq:FTc}.

Conversely, assume that $\bigl(T,xyz,(.\vert .)\bigr)$ is a
Freudenthal triple system. Then \eqref{eq:STSa} and
\eqref{eq:STSb} follow from \eqref{eq:FTa}, while \eqref{eq:FTa}
and \eqref{eq:FTb} imply \eqref{eq:STSd}. It is enough then to
prove that $d_{x,x}=[xx.]: y\mapsto xxy+\psi_{x,x}(y)$ is a
derivation of $(T,[...])$ for any $x\in T$, and since
$d_{T,T}\subseteq \frsp(T)$, it suffices to show that $d_{x,x}$ is
a derivation of the Freudenthal triple system or, equivalently,
that
\begin{equation}
d_{x,x}(yyz)=2(d_{x,x}(y))yz+yy(d_{x,x}(z))\label{eq:deri}
\end{equation}
for any $x,y,z\in T$. Linearizing \eqref{eq:FTc} in $y$, then
taking $z=y$ and using \eqref{eq:FTa} and \eqref{eq:FTb} one
obtains
\begin{multline}
2(xyt)xy+(yxx)yt+(xxt)yy+2(xyt\vert y)x\\
+(yxx\vert y)t +(txx\vert y)y+2(x\vert y)xyt+(t\vert
y)yxx=0.\label{eq:FTclin}
\end{multline}
Now, interchange $x$ and $y$ in \eqref{eq:FTclin} and subtract the
result from \eqref{eq:FTclin} to get
\begin{multline}
(xxy)yt-(yyx)xt+yy(xxt)-xx(yyt)-(x\vert yyt)x+(y\vert xxt)y\\
+4(x\vert y)xyt +(t\vert y)yxx-(t\vert x)xyy=0.\label{eq:FTcc}
\end{multline}
Adding \eqref{eq:FTcc} and \eqref{eq:FTc} with $z=t$ gives
\[
2(xxy)yt+yy(xxt)-xx(yyt)-2(x\vert yyt)x+4(x\vert y)xyt+2(x\vert
t)xyy=0,
\]
which is equivalent to \eqref{eq:deri}.
\end{proof}

\smallskip

\begin{remark}
There are several results in the literature relating different
triple systems. See, for instance, \cite[Section 2]{KamiyaIII},
\cite{FaulknerFerrar} or \cite{Brown}.

Related constructions of a $5$-graded Lie algebra from Freudenthal
triple systems (or equivalent ternary algebras) are given in
\cite{Faulkner,KanSkop}.
\end{remark}

\bigskip

The models of the exceptional Lie algebras in Section 3 fit in the
situation described in Theorem \ref{th:STS}. Therefore, by Theorem
\ref{th:FTS}, they provide models of  Freudenthal triple systems.

Take for instance our model for $\fre_8$ \eqref{eq:e8}. For any
$\sigma\in \calSe{8}\setminus\calSe{7}$, $8\in\sigma$ and
\eqref{eq:e8} can be rewritten as
\[
\fre_8=\Bigl( \fre_7\oplus \frsp(V_8)\Bigr)\bigoplus
\left(\Bigl(\bigoplus_{\sigma\in\calSe{8}\setminus\calSe{7}}
V(\sigma\setminus\{8\})\Bigr)\otimes V_8\right)\,,
\]
and we are in the situation of Theorem \ref{th:STS}. In
particular,
\[
\calT_{\fre_8,\fre_7}=
 \bigoplus_{\sigma\in\calSe{8}\setminus\calSe{7}} V(\sigma\setminus\{8\})
\]
is a simple Freudenthal triple system of dimension $2^3\times
7=56$, whose triple product and nondegenerate alternating bilinear
form can be easily computed in terms of the invariant maps
$\varphi_{\sigma,\tau}$ in \eqref{eq:phisigmatau}. In the same
vein,
\[
\begin{split}
\fre_7&=\Bigl( \frd_6\oplus \frsp(V_4)\Bigr)\bigoplus
\left(\Bigl(\bigoplus_{\sigma\in\calSe{7}\setminus\calS_{\frd_6}}
V(\sigma\setminus\{4\})\Bigr)\otimes V_4\right)\,,\\
\fre_6&=\Bigl( \fra_5\oplus \frsp(V_4)\Bigr)\bigoplus
\left(\Bigl(\bigoplus_{\sigma\in\calSf\setminus\calS_{\frc_3}}
\tilde V(\sigma\setminus\{4\})\Bigr)\otimes V_4\right)\,,\\
\frf_4&=\Bigl( \frc_3\oplus \frsp(V_4)\Bigr)\bigoplus
\left(\Bigl(\bigoplus_{\sigma\in\calSf\setminus\calS_{\frc_3}}
V(\sigma\setminus\{4\})\Bigr)\otimes V_4\right)\,,
\end{split}
\]
and we obtain simple Freudenthal triple systems defined on:
\[
\begin{split}
\calT_{\fre_7,\frd_6}&=
 \bigoplus_{\sigma\in\calSe{7}\setminus\calS_{\frd_6}} V(\sigma\setminus\{4\})\,,\\
\calT_{\fre_6,\fra_5}&=
 \bigoplus_{\sigma\in\calSf\setminus\calS_{\frc_3}} \tilde V(\sigma\setminus\{4\})\,,\\
\calT_{\frf_4,\frc_3}&=
 \bigoplus_{\sigma\in\calSf\setminus\calS_{\frc_3}} V(\sigma\setminus\{4\})\,,
\end{split}
\]
of dimension $2^3\times 4=32$, $2^3+2^2\times 3=20$ and
$2^3+2\times 3=14$, respectively.

The classification in \cite{Meyberg} and \cite{Brown} implies
that, over algebraically closed fields, these are the simple
Freudenthal triple systems associated to the split simple Jordan
algebras of degree $3$, that is, the Freudenthal triple systems
originally considered by Freudenthal.



\providecommand{\bysame}{\leavevmode\hbox
to3em{\hrulefill}\thinspace}
\providecommand{\MR}{\relax\ifhmode\unskip\space\fi MR }
\providecommand{\MRhref}[2]{%
  \href{http://www.ams.org/mathscinet-getitem?mr=#1}{#2}
} \providecommand{\href}[2]{#2}

\end{document}